\documentclass[12pt]{amsart}
\usepackage{amsmath,amsfonts,amscd,amssymb,graphics,psfrag}
\newtheorem{theorem}{Theorem}[section]
\newtheorem{proposition}[theorem]{Proposition}
\newtheorem{lemma}[theorem]{Lemma}

\theoremstyle{definition}

\theoremstyle{remark} \newtheorem{remark}[theorem]{Remark}
\numberwithin{equation}{section}
\newcommand{\field}[1]{\ensuremath{\mathbb{#1}}}
\newcommand{\CC}{\field{C}}
\newcommand{\DD}{\complex{D}}

\newcommand{\complex}[1]{\mathsf{#1}} 

\DeclareMathOperator{\dToda}{dToda}

\newcommand{\Z}{{\mathbb{Z}}}
\newcommand{\U}{\mathbb{H}}
\newcommand{\mL}{\mathcal{L}}

\newcommand{\C}{\mathbb{C}}
\newcommand{\Del}{\mathbb{D}}
\newcommand{\F}{\mathcal{F}}
\newcommand{\R}{\mathbb{R}}

\newcommand{\pa}{\partial}

\newcommand{\ov}{\overline}

\newcommand{\z}{\bar{z}}

\begin{document}
\title{L\"owner equations and dispersionless hierarchies}
\author{Takashi Takebe}
\address{Department of Mathematics, Ochanomizu University, Otsuka
2-1-1, Bunkyo-ku, Tokyo, 112-8610, Japan}
\email{takebe@math.ocha.ac.jp}
\author{Lee-Peng Teo}\address{Faculty of Information
Technology, Multimedia University, Jalan Multimedia, Cyberjaya,
63100, Selangor Darul Ehsan, Malaysia}\email{lpteo@mmu.edu.my}
\author{Anton Zabrodin}
\address{Institute of Biochemical Physics, Kosygina str. 4, 119991
Moscow, Russia and ITEP, Bol. Cheremushkinskaya str. 25, 117259
Moscow, Russia} \email{zabrodin@itep.ru}
\begin{abstract}
Using the Hirota representation of dispersionless dKP and dToda
hierarchies, we show that the chordal L\"owner equations and
radial L\"owner equations respectively serve as consistency
conditions for one variable reductions of these integrable
hierarchies. We also clarify the geometric meaning of this result
by relating it to the eigenvalue distribution of normal random
matrices in the large $N$ limit.
\end{abstract}
\maketitle
\section{Introduction}

The L\"owner equation is a differential equation obeyed by a
family of continuously varying univalent conformal maps
$G_{\lambda}(w)$
from the exterior of the unit circle
onto a domain with a slit
formed by a continuously increasing arc of a fixed curve.
The real parameter $\lambda$
characterizes the length of the arc.
Let us normalize the maps
so that $G_{\lambda}(w) =e^{\phi}w+ u_1 +
u_2 w^{-1}+\ldots$ as $w\to \infty$ with a real
$\phi =\phi (\lambda )$, then the L\"owner equation
reads
\begin{align}\label{intro1}
\frac{\pa G_{\lambda}(w)}{\pa \lambda}=-w\,
\frac{\pa G_{\lambda}(w)}{\pa w} \,\,
\frac{\sigma  (\lambda )+w}{\sigma  (\lambda )-w}\,\,
\frac{\pa \phi (\lambda )}{\pa \lambda}
\end{align}
The shape of the fixed
curve along which the endpoint of the varying slit runs is
encoded by the function $\sigma  (\lambda )$.
In fact $\sigma  (\lambda )$ is the pre-image
of the tip of the slit, so it lies on the unit circle:
$|\sigma  (\lambda )|=1$.
Equation (\ref{intro1}) was introduced by K.L\"owner
\cite{Lo} in 1923 in his attempt to prove the
Bieberbach conjecture \cite{Bie,deB,Duren}.
Now it is referred to as the \emph{radial L\"owner equation}.

An analog of equation (\ref{intro1}) called
\emph{chordal L\"owner equation} was introduced in 1999
\cite{ref10,Schramm} in other contexts:
\begin{align}\label{intro2}
\frac{\pa H_{\lambda}(w)}{\pa \lambda}=-\,\,
\frac{\pa H_{\lambda}(w)}{\pa w} \,\,
\frac{1}{U(\lambda )- w} \,\,
\frac{\pa a_1}{\pa \lambda}
\end{align}
Here $H_{\lambda}(w)=w+ a_1 w^{-1}+O(w^{-2})$ as $w\to \infty$
with a real coefficient $a_1$ and $U(\lambda )$ is a real-valued
function of $\lambda$. Under certain conditions, the function
$H_{\lambda}(w)$ conformally maps the upper half plane onto the
upper half plane with a slit, with $U(\lambda )$ being the
pre-image of the tip.

The Schramm's discovery \cite{Schramm,LSW} of the stochastic
L\"owner evolution (SLE) and its spectacular applications to the
conformal field theory
(see, e.g., the reviews \cite{Cardy,BB} and references therein)
have inspired a renewed interest in the theory of L\"owner equations.
One of its most interesting aspects is the relation to the integrable
hierarchies of nonlinear partial differential equations.
This relation was noticed by Gibbons and Tsarev \cite{ref10} yet before
the SLE boom. They have observed that the chordal L\"owner equation plays
a key role in classifying reductions of the dispersionless KP (dKP)
hierarchy \cite{TT2,TT1}. This point was further studied in
\cite{ref11,MMAM}. Recently, a similar role of the radial L\"owner equation
in the dispersionless coupled modified KP (dcmKP) hierarchy \cite{dcmKP}
was elucidated \cite{Manas,TTnew}. In the present paper we extend these results
to the dispersionless Toda (dToda) hierarchy \cite{TT4,TT1}.
Specifically, we characterize a class of reductions of the dToda hierarchy
in terms of solutions to the radial L\"owner equation. The representation of
the hierarchy in the Hirota form appears to be very useful for that purpose.
Using the Hirota framework, we also give a more transparent derivation of
the chordal L\"owner equation as a consistency condition for reductions of
the dKP hierarchy.

In short, our main result is the following\footnote{To make this
introduction as short as possible, we do not discuss the second
Lax function of the dToda hierarchy, for which similar results
hold true. For precise statements, see Propositions 5.6 and
5.11.}. Given any solution $G_{\lambda}(w)$ to the radial L\"owner
equation (\ref{intro1}) with any (continuous) $\sigma  (\lambda
)$, the Lax function $\mathcal{L}(w; \boldsymbol{t}
)=G_{\lambda}(w)$ solves the Lax equations of the dToda hierarchy
provided the dependence $\lambda =\lambda (\boldsymbol{t})$ is
given by the system of hydrodynamic type
\begin{align}\label{intro3}
\frac{\pa \lambda}{\pa t_n}=\xi_n (\lambda )\,
\frac{\pa \lambda}{\pa t_0}\,.
\end{align}
In these equations, the functions
$\xi_n (\lambda )$ are constructed in a canonical way
from $\mathcal{L}(w; \boldsymbol{t})$ and $\sigma  (\lambda )$.
In the context of the dKP and dcmKP hierarchies,
this result was earlier established in
\cite{ref10,MMAM,Manas,TTnew}.
Here we also prove the converse.
Let $\mathcal{L}(w; \boldsymbol{t})$ be the Lax function of the
dToda hierarchy.
Suppose one has a \emph{one-variable reduction} of the hierarchy,
i.e., the Lax function depends
on the hierarchical times $\boldsymbol{t}=
(\ldots, t_{-1}, t_0 , t_1 , \ldots )$
via a single function $\lambda (\boldsymbol{t})$:
$$
\mathcal{L}(w; \boldsymbol{t})=\mathcal{L}(w, \lambda (\boldsymbol{t}))\,.
$$
Then consistency of this ansatz
with the hierarchy implies that the function of two variables
$\mathcal{L}(w, \lambda )$ obeys the radial L\"owner equation
(\ref{intro1}) with some $\sigma  (\lambda )$ while the time dependence
of $\lambda$ is determined from system (\ref{intro3}).
A similar statement holds true for the dKP hierarchy
with the chordal L\"owner equation instead of the radial one.

The geometric meaning of this result and, more generally,
that of reductions of the dToda hierarchy, is clarified
in Section 6.
It is known
that solutions of the  dToda hierarchy obeying
certain reality conditions generate univalent conformal maps of
planar domains. The domain to be mapped (onto some fixed reference
domain) depends on which solution of the hierarchy one takes as
well as on values of the hierarchical ``times'' $t_n$. Generic
solutions correspond to mappings of planar domains bounded by
smooth simple curves, with the times being suitably defined
moments of the domain \cite{ref6_1,ref6_2,ref7,ref8,ref9}.
Non-generic solutions, or {\it reductions} of the hierarchy, are
known to be related to conformal maps of domains with highly
singular boundaries, such as slit domains \cite{ref10,ref11}.

Both types of solutions and conformal maps corresponding to them
can be understood in terms of the model of normal random matrices
\cite{CY,ref20} in the large $N$ limit. Actually, we use only the
eigenvalue integral representation, which is the partition
function of the 2D Coulomb gas in external field. The matrix
models are known to be solvable (even at finite $N$) in terms of
integrable hierarchies, which become dispersionless as $N\to
\infty$. At the same time, the domains subject to conformal maps
under study appear in this context as complements of the supports
of eigenvalues in the $N\to \infty$ matrix model. This picture
provides a transparent geometric meaning of the reductions and
conformal maps of slit domains corresponding to them.

\section{L\"owner Equations}
\subsection{Radial L\"owner Equation}
The L\"owner theory was introduced by L\"owner in \cite{Lo}. It
plays an important role in solving the Bieberbach conjecture
\cite{Bie, deB, Wein}. Here we briefly recall the necessary
ingredients  in the context we need.  One can refer to \cite{Pom,
Duren, Conway} for a complete exploration of the theory.

\begin{figure}[h]
 \begin{center}
  \psfrag{Be}{$\sf{B}$}
  \psfrag{0}{$0$}
  \psfrag{Ga}{$\Gamma$}
  \psfrag{Gal}{$\Gamma_\lambda$}
  \psfrag{Fl}{$z=F_\lambda(w)$}
  \psfrag{fl}{$f_\lambda(z)=w$}
  \psfrag{e(l)}{$\eta(\lambda)$}
  \psfrag{Del}{$\Del$}
  \includegraphics{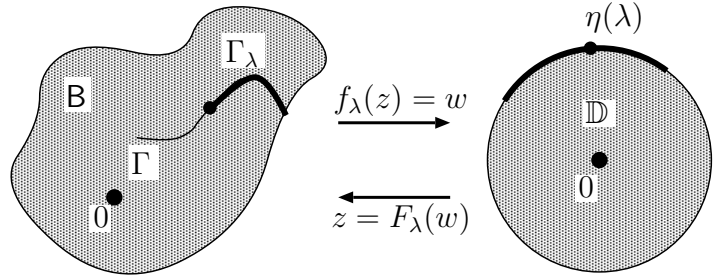}
 \end{center}
\label{fig:2-1-1}
\caption{Conformal maps $F_\lambda$, $f_\lambda$.}
\end{figure}

Let $\Del$ be the unit disc and $\Del^*$ its exterior. Let
$\sf{B}$ be a simply connected domain in the complex plane $\C$
containing the origin and let $\Gamma:[0,\infty)\rightarrow \C$ be
a  Jordan arc growing inside $\sf{B}$, i.e., $\Gamma(0)$ is a
point lying on the boundary of $\sf{B}$, and for $\mu\in
(0,\infty)$, $\Gamma(\mu)\in \sf{B}$. Let
$\Gamma_{\lambda}=\left.\Gamma\right|_{[0,\lambda]}$ be the
arc of $\Gamma$ between $0$ and $\lambda$ and let
$\sf{B}_{\lambda}=\sf{B}\setminus$   $\Gamma_{\lambda}$ be the
corresponding slit domain. By the Riemann mapping theorem, there
exists a unique conformal map  $F_{\lambda}$  mapping the unit
disc $\Del$ onto the domain $\sf{B}_{\lambda}$ and satisfying
$F_{\lambda}(0)=0$, $F_{\lambda}'(0)=e^{-\phi(\lambda)}>0$.
Let $\eta(\lambda)$ be the point on $S^1 =\partial \Del$
mapped to the tip of
$\Gamma_{\lambda}$ by $F_{\lambda}$, i.e.,
$F_{\lambda}(\eta(\lambda))=\Gamma(\lambda)$. Then
$F(w,\lambda)=F_{\lambda}(w)$ satisfies the following differential
equation:
\begin{align}\label{radialin}
\frac{\pa F(w,\lambda)}{\pa \lambda}=-w\frac{\pa F(w,\lambda)}{\pa
w}\frac{\eta(\lambda)+w}{\eta(\lambda)-w}\frac{\pa\phi(\lambda)}{\pa\lambda},
\hspace{1cm}w\in\Del, \lambda\in[0,\infty),
\end{align}
called the (radial) L\"owner equation. Let $f_{\lambda}$ be the
inverse function of $F_{\lambda}$. The corresponding L\"owner
equation for $f(z,\lambda)=f_{\lambda}(z)$ is
\begin{align*}
\frac{\pa f(z,\lambda)}{\pa \lambda}=
f(z,\lambda)\frac{\eta(\lambda)+f(z,\lambda)}{\eta(\lambda)-
f(z,\lambda)}\frac{\pa\phi(\lambda)}{\pa\lambda}
,\hspace{1cm}z\in \sf{B}_{\lambda}, \lambda\in[0,\infty).
\end{align*}
(See Figure~\ref{fig:2-1-1}.)

Conversely, given  a measurable
function $\eta:[0,\infty)\rightarrow S^1$
and an increasing
differentiable function $\phi:[0,\infty)\rightarrow \R$, 
the L\"owner equation \eqref{radialin} with the initial
condition $F_0(\Del)=\sf{B}$,
$F_0(0)=0$, $F_{0}'(0)=e^{-\phi(0)}$ always has a unique solution
$F_{\lambda}:\Del\rightarrow \C$ which maps the unit
disc to a family of
continuously shrinking domains
$F_{\lambda}(\Del)$\footnote{However, $F_{\lambda}(\Del)$
are not necessarily slit domains, see \cite{Kufarev,MarshallRohde}.}.

\begin{figure}[h]
 \begin{center}
  \psfrag{Be}{$\sf{B}$}
  \psfrag{0}{$0$}
  \psfrag{oo}{$\infty$}
  \psfrag{Ga}{$\Gamma$}
  \psfrag{Gal}{$\Gamma_\lambda$}
  \psfrag{Gl}{$z=G_\lambda(w)$}
  \psfrag{gl}{$g_\lambda(z)=w$}
  \psfrag{s(l)}{$\sigma(\lambda)$}
  \psfrag{Del}{$\Del$}
  \includegraphics{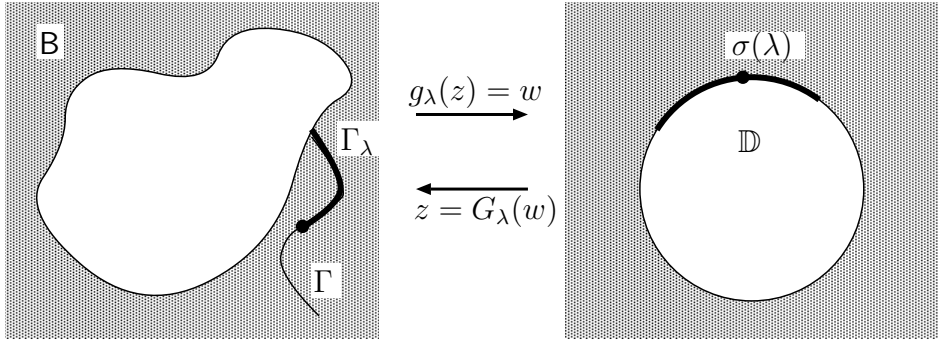}
 \end{center}
\label{fig:2-1-2}
\caption{Conformal maps $G_\lambda$, $g_\lambda$.}
\end{figure}

In case when $\sf{B}$ is a simply connected domain in $\hat{\C}$ that
contains $\infty$ and $\Gamma:[0,\lambda)\rightarrow \hat{\C}$ is a
growing Jordan arc in $\sf{B}$, let
$\Gamma_{\lambda}=\left.\Gamma\right|_{[0,\lambda]}$, and let
$G_{\lambda}$ be the unique conformal map from the exterior disc
$\Del^*$ onto $\sf{B}_{\lambda}=\sf{B}\setminus \Gamma_{\lambda}$
normalized such that $G_{\lambda}(\infty)=\infty$,
$G_{\lambda}'(\infty)=e^{\phi(\lambda)}>0$ and let $g_{\lambda}$ be the
inverse function of $G_{\lambda}$. Then $G(w,\lambda)=G_{\lambda}(w)$
and $g(z,\lambda)=g_{\lambda}(z)$ satisfy the L\"owner equations
\begin{align}
\frac{\pa G(w,\lambda)}{\pa \lambda}=&-w\frac{\pa
G(w,\lambda)}{\pa
w}\frac{\sigma(\lambda)+w}{\sigma(\lambda)-w}\frac{\pa\phi(\lambda)}{\pa\lambda}
, \hspace{1cm}w\in\Del^*,\lambda\in[0,\infty),\label{radialout}\\
\frac{\pa g(z,\lambda)}{\pa \lambda}=&
g(z,\lambda)\frac{\sigma(\lambda)+
g(z,\lambda)}{\sigma(\lambda)-g(z,\lambda)}\frac{\pa\phi(\lambda)}{\pa\lambda}
,\hspace{1.5cm}z\in \sf{B}_{\lambda},
\lambda\in[0,\infty).\nonumber
\end{align}
(See Figure~\ref{fig:2-1-2})
Here $\sigma(\lambda)$ is the point on $S^1$ mapped to
the tip of $\Gamma_{\lambda}$ by $G_{\lambda}$.

In the special case when the domain $\sf{B}$ is the exterior disc
$\Del^*$, let $J_{\lambda}\subset
S^1$ be the pre-image of
$\Gamma_{\lambda}$ under the map $G_{\lambda}$.
Then $G_{\lambda}$ maps the arc $S^1\setminus J_{\lambda}$
to an arc of $S^1$. By the Schwarz reflection principle, the map $G_{\lambda}$
has an analytic continuation $\tilde{G}_{\lambda}$ defined on
$\Del\cup \Del^* \cup (S^1\setminus J_{\lambda})$, with the image
$\hat{\C}\setminus (\Gamma_{\lambda}\cup
\tilde{\Gamma}_{\lambda})$, where $\tilde{\Gamma}_{\lambda}$ is
the mirror image of $\Gamma_{\lambda}$ under reflection
in the unit circle, i.e.
$\tilde{\Gamma}_{\lambda}=\{1/\z\;:\; z\in\Gamma_{\lambda}\}$. The
map $F_{\lambda}=\left.\tilde{G}_{\lambda}\right|_{\Del}$ is given
explicitly by
$$F_{\lambda}(w) =\ov{G_{\lambda}(1/\bar{w})}^{\;-1}.$$
Using this relation, it is easy
to check that $F(w,\lambda)=F_{\lambda}(w)$ satisfies the L\"owner
equation \eqref{radialin} with $\eta(\lambda)=\sigma(\lambda)$.
Therefore, the differential equations for $G(w,\lambda)$ and
$F(w,\lambda)$ are the same but defined in different domains.

\subsection{Chordal L\"owner Equation} 
\begin{figure}[h]
 \begin{center}
  \psfrag{Ga}{$\Gamma$}
  \psfrag{Gal}{$\Gamma_\lambda$}
  \psfrag{Hl}{$z=H_\lambda(w)$}
  \psfrag{hl}{$h_\lambda(z)=w$}
  \psfrag{U(l)}{$U(\lambda)$}
  \psfrag{H}{$\U$}
  \includegraphics{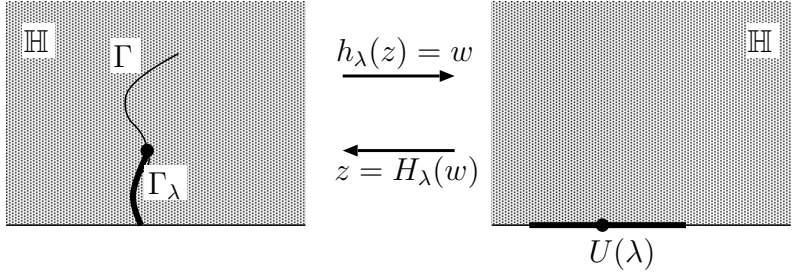}
 \end{center}
\label{fig:2-2}
\caption{Conformal maps $H_\lambda$, $h_\lambda$.}
\end{figure}

The chordal L\"owner equation
\cite{Schramm,LSW,ref10} is an
analogue of the radial L\"owner equation for conformal maps on the
upper half plane $\U$. We say that a conformal map
$H:\U\rightarrow \U$ that has continuous extension to $\pa\U$ is
normalized with respect to the point $\infty$ if it has a Laurent
series expansion of the form
$$H(w)=w+\frac{a_1}{w}+\frac{a_2}{w^2}+\ldots.$$
Let $H_\lambda:\U\rightarrow \U$ be a sequence
of normalized conformal maps such that $H_\lambda(\U)=\U\setminus
\Gamma_\lambda$, where $\Gamma: [0,\infty)\rightarrow \C$ is a
growing Jordan arc in $\U$ and
$\Gamma_{\lambda}=\left.\Gamma\right|_{[0,\lambda]}$.
 Denote by
$U(\lambda)$ the point on $\R$ that is mapped to the tip of
$\Gamma_{\lambda}$ by $H_{\lambda}$ and let $h_{\lambda}$ be the
inverse function of $H_{\lambda}$. Then
$H(w,\lambda)=H_\lambda(w)$ and $h(z,\lambda)=h_{\lambda}(z)$
satisfy the differential equations:
\begin{align}\label{chordal}
\frac{\pa H(w,\lambda)}{\pa\lambda}=-\frac{\pa H(w,\lambda)}{\pa
w}\frac{1}{U(\lambda)-w}\frac{\pa
a_1(\lambda)}{\pa \lambda},\hspace{1cm} \lambda\in [0,\infty),\\
\frac{\pa h(z,\lambda)}{\pa
\lambda}=\frac{1}{U(\lambda)-h(z,\lambda)}\frac{\pa
a_1(\lambda)}{\pa \lambda},\hspace{1cm}\lambda\in
[0,\infty),\nonumber
\end{align}
which are called the chordal L\"owner equations for $H_{\lambda}$
and $h_{\lambda}$ respectively. (See Figure~\ref{fig:2-2}.)

\section{ Grunsky coefficients and Faber polynomials  }
In this section we review the definitions of Faber polynomials and
Grunsky coefficients. For details, see \cite{Pom, Duren, Teo}. Let
\begin{align*} f(z)=& rz+\alpha_2 z^2+\alpha_3
z^3+\ldots\;\;\text{and}\hspace{0.5cm} g(z)=r^{-1} z + \beta_0 +
\frac{\beta_1}{z}+\frac{\beta_2}{z^2}+\ldots
\end{align*}be functions univalent in a neighbourhood of $0$ and
 $\infty$ respectively. The  Grunsky coefficients $b_{m,n}$, $m,n\in \Z$, of
 the pair $(f,g)$ are defined by the expansions
\begin{align}\label{Grunsky}
\log\frac{g(z_1)-g(z_2)}{z_1-z_2}=& -\log
r-\sum_{m=1}^{\infty}\sum_{n=1}^{\infty} b_{m,n} z_1^{-m}z_2^{-n}\\
\log \frac{g(z_1)-f(z_2)}{z_1}=&-\log
r-\sum_{m=1}^{\infty}\sum_{n=0}^{\infty} b_{m,-n}z_1^{-m}z_2^n\nonumber\\
\log\frac{f(z_1)-f(z_2)}{z_1-z_2}=&-
\sum_{m=0}^{\infty}\sum_{n=0}^{\infty}b_{-m,-n}z_1^mz_2^n\nonumber
\end{align}and $b_{-m,n}=b_{n,-m}$ for $n\geq 1, m\geq 0$. By
definition, $b_{m,n}=b_{n,m}$ for all $n,m\in \Z$, $b_{0,0}=-\log
r$ and
\begin{align}\label{Grunsky3}
\log\frac{g(z)}{z}=-\log r-\sum_{n=1}^{\infty}
b_{n,0}z^{-n},\hspace{1cm}\log\frac{f(z)}{z}=\log r
-\sum_{n=1}^{\infty}b_{-n,0}z^n.
\end{align}Using these formulas, the second and the third equations in
\eqref{Grunsky} can be rewritten as \begin{align}\label{Grunsky2}
\log \left(1-\frac{f(z_2)}{g(z_1)}\right)
=&-\sum_{m=1}^{\infty}\sum_{n=1}^{\infty} b_{m,-n}z_1^{-m}z_2^n\\
\log\frac{f(z_1)^{-1}-f(z_2)^{-1}}{z_1^{-1}-z_2^{-1}}=&-\log
r-\sum_{m=1}^{\infty}\sum_{n=1}^{\infty}b_{-m,-n}z_1^mz_2^n\nonumber.
\end{align}
 For $g$ alone,
$b_{m,n}$, $m,n\geq 0$ are called the Grunsky coefficients of $g$.

The Faber polynomials $\Phi_n(w), n\geq 1$ of $g$ and Faber
polynomials $\Psi_n, n\geq 1$ of $f$ are defined by
\begin{align}\label{Faber}
\log\frac{g(z)-w}{r^{-1}z}=&-\sum_{n=1}^{\infty}
\frac{\Phi_n(w)}{n}z^{-n},\\\log\frac{w-f(z)}{w}=&
\log
\frac{f(z)}{rz}-\sum_{n=1}^{\infty}\frac{\Psi_n(w)}{n}z^n,\nonumber
\end{align}
so that
\begin{align*}
\Phi_n(w)=(G(w)^n)_{\geq
0},\hspace{1cm}\Psi_n(w)=(F(w)^{-n})_{\leq 0},
\end{align*}
where $F$ and $G$ are the inverse functions of $f$ and $g$
respectively, and for
$S\subset \Z$ we denote $(\sum_{k\in \Z} A_k w^k)_{S}=\sum_{k\in
S}A_kw^k$. Taking derivatives with respect to  $\log w$, we have
\begin{align}\label{Fa1}
\frac{w}{w-g(z)}=-\sum_{n=1}^{\infty}\frac{w\Phi_n'(w)}{n}z^{-n},
\hspace{1cm}\frac{f(z)}{w-f(z)}
=-\sum_{n=1}^{\infty}\frac{w\Psi_n'(w)}{n}z^n.
\end{align}

\section{Dispersionless hierarchies}
We briefly review the dispersionless Toda (dToda) hierarchy,
the dispersionless
KP (dKP) hierarchy and their tau functions. For details, see
\cite{TT2, TT4, TT1}.

\subsection{Dispersionless KP hierarchy}
The fundamental quantity in the dKP hierarchy is a formal power series
\begin{align}\label{series2}
\mL(w;\boldsymbol{t})= w+ \sum_{n=1}^{\infty}
u_{n+1}(\boldsymbol{t})w^{-n}
\end{align}with coefficients depending on the independent
variables $\boldsymbol{t}=(x+t_1, t_2, t_3,\ldots)$. The Lax
representation of the dKP hierarchy is
\begin{align}\label{Lax2}
\frac{\pa \mL}{\pa t_n}=\{\mathcal{B}_n,\mL\}, \hspace{1cm}
\mathcal{B}_n=(\mL^n)_{\geq 0}.
\end{align}Here $\{\,\cdot\,,\,\cdot\,\}$ is the Poisson bracket
$$\{h_1, h_2\}=\frac{\pa h_1}{\pa w}\frac{\pa h_2}{\pa
x}-\frac{\pa h_1}{\pa x}\frac{\pa h_2}{\pa w}.$$Let
$k(z;\boldsymbol{t})$ be the inverse of $\mL(w;\boldsymbol{t})$,
i.e. $k(\mL(w;\boldsymbol{t});\boldsymbol{t})=w$ and
$\mL(k(z;\boldsymbol{t});\boldsymbol{t})=z$. There exists a tau
function $\tau_{\text{dKP}}(\boldsymbol{t})$ that generates the
Grunsky coefficients $b_{m,n}(\boldsymbol{t})$, $m,n\geq 1$ of
$k(z;\boldsymbol{t})$. More precisely,
\begin{align}\label{identity2}
\frac{\pa^2\log\tau_{\text{dKP}}(\boldsymbol{t})}{\pa t_m\pa t_n}
=-mnb_{m,n}(\boldsymbol{t}).\end{align}
Introduce the operator
$$D(z)=\sum_{n=1}^{\infty}\frac{z^{-n}}{n}\frac{\pa}{\pa t_n},$$
and set $\F=\log\tau_{dKP}$
(the ``free energy''). Then the first equation of \eqref{Grunsky} can be
written as
\begin{align}\label{equation2}
\log\frac{k(z_1)-k(z_2)}{z_1-z_2}=&D(z_1)D(z_2)\F,\hspace{1cm}
k(z)= z-D(z)\pa_{t_1}\F,
\end{align}
which is the Hirota equation for the dKP hierarchy \cite{TT1,
ref15, ref13_1,ref13_2} . Conversely, it was proved  \cite{ref1,
ref13_1,ref13_2,Teo} that the Hirota equation implies the dKP
hierarchy. We reformulate this statement below in the form we need
later.
\begin{proposition}\label{tau1}
Let $k(z;\boldsymbol{t})$  be the inverse of the formal power
series $\mL(w;\boldsymbol{t})$ of the form \eqref{series2} and let
$b_{m,n}(\boldsymbol{t})$ be the Grunsky coefficients of
$k(z;\boldsymbol{t})$. If there exists a function
$\tau_{\text{dKP}}(\boldsymbol{t})$ satisfying \eqref{identity2},
then $\mL(w;\boldsymbol{t})$ is a solution of the dKP hierarchy
\eqref{Lax2}.
\end{proposition}

\subsection{Dispersionless Toda hierarchy}

The fundamental quantities in the dToda hierarchy are two formal power
series in $w$\;:
\begin{align}\label{series1}
\mL(w;\boldsymbol{t})=& r(\boldsymbol{t})w + \sum_{n=0}^{\infty}
u_{n+1}(\boldsymbol{t}) w^{-n},\hspace{0.3cm}
\tilde{\mL}^{-1}(w;\boldsymbol{t}) =& r(\boldsymbol{t})w^{-1} +
\sum_{n=0}^{\infty} \tilde{u}_{n+1}(\boldsymbol{t}) w^n.
\end{align}
Here $u_n(\boldsymbol{t})$ and $\tilde{u}_n(\boldsymbol{t})$ are
functions of the independent variables $t_n, n \in \Z$, which we
denote collectively by $\boldsymbol{t}$. The Lax representation is
\begin{align}\label{Lax1}
\frac{\pa \mL}{\pa t_n} =\{ \mathcal{B}_n, \mL\}_T, \hspace{2cm}
\frac{\pa \mL}{\pa t_{-n}} = \{\tilde{\mathcal{B}}_n, \mL\}_T,\\
\frac{\pa \tilde{\mL}}{\pa t_n} =\{ \mathcal{B}_n,
\tilde{\mL}\}_T, \hspace{2cm} \frac{\pa \tilde{\mL}}{\pa t_{-n}} =
\{\tilde{\mathcal{B}}_n, \tilde{\mL}\}_T.\nonumber
\end{align}
Here \begin{align*}
\mathcal{B}_n=(\mL^n)_{>0}+\frac{1}{2}(\mL^n)_0,\hspace{1cm}
\tilde{\mathcal{B}}_n=
(\tilde{\mL}^{-n})_{<0}+\frac{1}{2}(\tilde{\mL}^{-n})_{0},
\end{align*} and $\{ \cdot, \cdot\}_T$ is the Poisson bracket for the dToda
hierarchy
\[
\{h_1, h_2\}_T =w \frac{\pa h_1}{\pa w} \frac{\pa h_2}{\pa
t_0}-w\frac{\pa h_1}{\pa t_0} \frac{\pa h_2}{\pa w}.
\]

Let  $p(z;\boldsymbol{t})$ and $\tilde{p}(z;\boldsymbol{t})$  be
the inverses of $\mL(w;\boldsymbol{t})$ and
$\tilde{\mL}(w;\boldsymbol{t})$ respectively, i.e.
$p(\mL(w;\boldsymbol{t});\boldsymbol{t})=w,
\mL(p(z;\boldsymbol{t});\boldsymbol{t})=z,
\tilde{p}(\tilde{\mL}(w;\boldsymbol{t});\boldsymbol{t})=w$ and
$\tilde{\mL}(\tilde{p}(z;\boldsymbol{t});\boldsymbol{t})=z$. There
exists a tau function $\tau_{\dToda}(\boldsymbol{t})$ which
generates the Grunsky coefficients $b_{m,n}(\boldsymbol{t})$,
$m,n\in \Z$ of $(\tilde{p}(z;\boldsymbol{t}),
p(z;\boldsymbol{t}))$. More precisely,
\begin{align}\label{identity1}
\frac{\pa^2\log\tau_{\dToda}(\boldsymbol{t})}{\pa t_m\pa
t_n}=\begin{cases}-|mn|b_{mn}(\boldsymbol{t})
,\hspace{1cm}&\text{if}\;\;m\neq 0, n\neq 0\\
|m|b_{m,0}(\boldsymbol{t}),&\text{if}\;\;m\neq 0, n=0\\
-2b_{00}(\boldsymbol{t}),&\text{if}\;\;m=n=0.
\end{cases}
\end{align}
Let us introduce the operators
\begin{align}\label{op}D(z)=\sum_{n=1}^{\infty}\frac{z^{-n}}{n}
\frac{\pa}{\pa t_n},
\hspace{1cm}\tilde{D}(z)=\sum_{n=1}^{\infty}\frac{z^n}{n}\frac{\pa}{\pa
t_{-n}}\end{align}
and  define the free energy $\F$ by
$\F=\log\tau_{\dToda}$. Then we can rewrite the first equation in
\eqref{Grunsky}, equations \eqref{Grunsky2} and \eqref{Grunsky3}
in the form
\begin{align}\label{equation1}
&\log \frac{p(z_1)-p(z_2)}{z_1-z_2}=-\frac{1}{2}
\frac{\pa^2\F}{\pa t_0^2}+D(z_1)D(z_2)\F\\
&\log\left(1-\frac{\tilde{p}(z_2)}{p(z_1)}\right)=
D(z_1)\tilde{D}(z_2)\F\nonumber\\
&\log\frac{\tilde{p}(z_1)^{-1}-\tilde{p}(z_2)^{-1}}{z_1^{-1}-z_2^{-1}}=
-\frac{1}{2}\frac{\pa^2\F}{\pa t_0^2}+
\tilde{D}(z_1)\tilde{D}(z_2)\F\nonumber\\
&\log \frac{p(z)}{z}=-\frac{1}{2}\frac{\pa^2\F}{\pa
t_0^2}-D(z)\pa_{t_0}\F,\hspace{1cm}\log\frac{\tilde{p}(z)}{z}=
\frac{1}{2}\frac{\pa^2\F}{\pa
t_0^2}-\tilde{D}(z)\pa_{t_0}\F.\nonumber
\end{align}
This is the system of Hirota
equations for the dToda hierarchy
\cite{ref6_1,ref6_2,ref8,ref9,ref1}.
Conversely, it was shown in \cite{ref1,ref13_1,ref13_2,Teo}
that the Hirota equations imply the dToda
hierarchy. We reformulate this statement below in the form we need.
\begin{proposition}\label{tau2}
Let $p(z;\boldsymbol{t})$ and $\tilde{p}(z;\boldsymbol{t})$ be the
inverses of the formal power series $\mL(w;\boldsymbol{t})$ and
$\tilde{\mL}(w;\boldsymbol{t})$ of the form \eqref{series1} and
let $b_{m,n}(\boldsymbol{t})$ be the Grunsky coefficients of
$(\tilde{p}(z;\boldsymbol{t}),p(z;\boldsymbol{t}))$. If there
exists a function $\tau_{\dToda}(\boldsymbol{t})$ satisfying
\eqref{identity1}, then $(\mL(w;\boldsymbol{t}),
\tilde{\mL}(w;\boldsymbol{t}))$ is a solution of the dToda
hierarchy \eqref{Lax1}.
\end{proposition}

\section{One variable reductions of dispersionless hierarchies }

In this section, we consider  the one-variable reductions of the
dKP and dToda hierarchies from the perspective of their Hirota
equations.

\subsection{One variable reduction of dKP hierarchy}

\begin{proposition}\label{Pro1}
Suppose that $\mL(w;\boldsymbol{t})$
is a solution of the dKP
hierarchy whose dependence on $\boldsymbol{t}=(x+t_1, t_2,
t_3,\ldots)$ is through a single variable $\lambda$. Namely, there
exists a function $\lambda(\boldsymbol{t})$ of $\boldsymbol{t}$ such
that
$$\mL(w;\boldsymbol{t})=\mL(w,\lambda(\boldsymbol{t}))=w+\sum_{n=1}^{\infty}
u_{n+1}(\lambda(\boldsymbol{t}))w^{-n}.$$
Then $\mL(w,\lambda)$
satisfies the chordal L\"owner equation \eqref{chordal} with
respect to $\lambda$. Moreover, let $k(z,\lambda)$ be the
inverse function of $\mL(w,\lambda)$, $\Phi_n(w,\lambda)$, $n\geq
1$ be the Faber polynomials of $k(z,\lambda)$ and
$\chi_n(\lambda)=\Phi_n'(U(\lambda),\lambda)$, then
$\lambda(\boldsymbol{t})$ satisfies the hydrodynamic type equation
\begin{align}\label{hodograph}
\frac{\pa\lambda}{\pa t_n}=\chi_n(\lambda)\frac{\pa \lambda}{\pa
t_1}.
\end{align}

\end{proposition}

\begin{proof}
 Let $\F(\boldsymbol{t})=\log\tau_{\text{dKP}}(\boldsymbol{t})$
be the free energy of the solution
$\mL(w,\lambda(\boldsymbol{t}))$. From the second equation in
\eqref{equation2} we have
\begin{align*}
D(z_1)k(z_2,\lambda)=D(z_2)k(z_1,\lambda).
\end{align*}
Setting $z_1=z$ and comparing the coefficient of $z_2^{-1}$ on both
sides, we obtain:
\begin{align*}
D(z)u =
-\pa_{t_1}k(z,\lambda)=-\pa_{\lambda}k(z,\lambda)\pa_{t_1}\lambda,
\hspace{0.6cm}\text{where}\;\;\;u=\frac{\pa^2\F}{\pa
t_1^2}=u_2(\lambda(\boldsymbol{t})).
\end{align*}
Suppose $\pa_{\lambda} u\neq 0$, then we have by the chain rule:
\begin{align}\label{re1}
 D(z)\lambda=-\pa_{\lambda}k(z,\lambda)\frac{\pa_{t_1}\lambda}{\pa_{\lambda}
 u}\end{align}and\begin{align*}D(z_2)k(z_1,\lambda)=-\pa_{\lambda}k(z_1,\lambda)
 \pa_{\lambda}k(z_2,\lambda)
\frac{\pa_{t_1}\lambda}{\pa_{\lambda} u}.
\end{align*}
Applying $\pa_{t_1}$ to both sides of the first equation in
\eqref{equation2}, we get
\begin{align*}
\frac{\pa_{\lambda}k(z_1,\lambda)-
\pa_{\lambda}k(z_2,\lambda)}{k(z_1,\lambda)-k(z_2,\lambda)}\pa_{t_1}\lambda
=&D(z_2)(z_1-k(z_1,\lambda))=\pa_{\lambda}k(z_1,\lambda)
\pa_{\lambda}k(z_2,\lambda)
\frac{\pa_{t_1}\lambda}{\pa_{\lambda} u}.
\end{align*}If $\pa_{t_1}\lambda\neq 0$, then
\begin{align*}
\left(\pa_{\lambda}k(z_1,\lambda)-\pa_{\lambda}k(z_2,\lambda)
\right)\pa_{\lambda}u=\pa_{\lambda}k(z_1,\lambda)
\pa_{\lambda}k(z_2,\lambda)(k(z_1,\lambda)-k(z_2,\lambda)),
\end{align*}
which means that
\begin{align}
\label{defineU}(\pa_{\lambda}k(z,\lambda))^{-1}\pa_{\lambda}
u+k(z,\lambda)=U(\lambda)\end{align} is independent of $z$. This
gives
$$\frac{\pa k(z,\lambda)}{\pa
\lambda}=\frac{1}{U(\lambda)-k(z,\lambda)}\frac{\pa
u}{\pa\lambda},$$ which is the chordal L\"owner equation
\eqref{chordal} for $k(z,\lambda)$.

Now,  from equation \eqref{re1} and the first equation in
\eqref{Fa1}, we have
\begin{align*}
\sum_{n=1}^{\infty} \frac{z^{-n}}{n}\frac{\pa\lambda}{\pa
t_n}=\frac{1}{k(z,\lambda)-U(\lambda)}\frac{\pa\lambda}{\pa
t_1}=\sum_{n=1}^{\infty}\frac{z^{-n}}{n}\Phi_n'(U(\lambda),
\lambda)\frac{\pa\lambda}{\pa
t_1}.
\end{align*} Therefore, with
$\chi_n(\lambda)=\Phi_n'(U(\lambda),\lambda)$, we have
\begin{align*}\frac{\pa\lambda}{\pa t_n}=\chi_n(\lambda)\frac{\pa\lambda}{\pa
t_1}.\end{align*}\end{proof}
\begin{remark}
If we assume that the function $\mL(w;\boldsymbol{t})$ maps
 the upper half plane $\U$ onto a slit sub-domain of the upper half plane,
then by
 letting $z\rightarrow \pa \U$ in \eqref{defineU}, we find that
 $U(\lambda)\in \R$.
\end{remark}
Proposition \ref{Pro1} implies that for
$\mL(w,\lambda(\boldsymbol{t}))$  to satisfy the dKP hierarchy,
$\lambda(\boldsymbol{t})$ must satisfy the hydrodynamic type
equation \eqref{hodograph}. In fact, \eqref{hodograph}  can be
solved by the general hodograph method of Tsarev \cite{Tsarev}:
\begin{lemma}\label{le1}
Let $R(\lambda)$ be any function of $\lambda$. If
$\lambda(\boldsymbol{t})$ is defined implicitly by the hodograph
relation
\begin{align}\label{ho2}
x+t_1+\sum_{n=2}^{\infty}\chi_n(\lambda)t_n =R(\lambda),
\end{align}
then $\lambda(\boldsymbol{t})$ satisfies \eqref{hodograph}.
\end{lemma}
\begin{proof}
A straightforward computation.
\end{proof}
\begin{remark}
If all the coefficients of the series $\mL(w,\lambda)$ are real,
then $\Phi_n(w,\lambda)$ is a polynomial in $w$ with real
coefficients. Therefore, imposing the condition that all the
variables $t_n$ are real, the relation \eqref{ho2} can be solved
for $\lambda$ as a real-valued function in a certain domain of
$\boldsymbol{t}$.
\end{remark}
As a  converse to Proposition \ref{Pro1}, Gibbons and Tsarev
\cite{ref11}, Yu and Gibbons \cite{ref10}, Man{\~n}as,
Mart{\'{\i}}nez~Alonso and Medina \cite{MMAM} and others
have shown that a solution of the chordal L\"owner equation
\eqref{chordal} together with equation \eqref{hodograph} give rise
to a solution of the dKP hierarchy. Here we give an independent proof
using the Hirota equations.

\begin{proposition}
Let $H(w,\lambda)$ be a solution of the chordal L\"owner equation
\eqref{chordal} and $\lambda(\boldsymbol{t})$ a solution of
equation \eqref{hodograph}. Then
$\mL(w;\boldsymbol{t})=H(w,\lambda(\boldsymbol{t}))$ is a solution
of the dKP hierarchy \eqref{Lax2}.
\end{proposition}\begin{proof}
Let $h(z,\lambda)$ be the inverse function of $H(w,\lambda)$ and
let $b_{m,n}(\boldsymbol{t})=b_{m,n}(\lambda(\boldsymbol{t}))$,
$m,n\geq 1$ be the Grunsky coefficients of $h(z,\lambda)$.
Applying $\pa_{t_k}$ to both sides of the first equation of
\eqref{Grunsky} and using the L\"owner equation for $h(z,\lambda)$
\eqref{chordal}, equation \eqref{hodograph} and the first equation
in \eqref{Fa1}, we have
\begin{align*}
&-\sum_{m=1}^{\infty}\sum_{n=1}^{\infty}\frac{\pa
b_{m,n}(\boldsymbol{t})}{\pa t_k}z_1^{-m}z_2^{-n}\\=&
\frac{\pa_{\lambda}h(z_1,\lambda)-
\pa_{\lambda}h(z_2,\lambda)}{h(z_1,\lambda)-h(z_2,\lambda)}
\frac{\pa\lambda}{\pa t_k}
\\=&\frac{1}{(U(\lambda)-h(z_1))(U(\lambda)-h(z_2))}\frac{\pa
a_1(\lambda)}{\pa\lambda}\chi_k(\lambda)\frac{\pa\lambda}{\pa
t_1}\\
=&\sum_{m=1}^{\infty}\sum_{n=1}^{\infty}
\frac{1}{mn}\chi_m(\lambda)\chi_n(\lambda)z_1^{-m}z_2^{-n}\frac{\pa
a_1(\lambda)}{\pa\lambda}\chi_k(\lambda)\frac{\pa\lambda}{\pa
t_1}.
\end{align*}
Comparing coefficients gives
\begin{align*}
-mn\frac{\pa b_{m,n}(\boldsymbol{t})}{\pa t_k}
=\chi_m(\lambda)\chi_n(\lambda)\chi_k(\lambda)\frac{\pa
a_1(\lambda)}{\pa\lambda}\frac{\pa\lambda}{\pa t_1},
\end{align*}
which is completely symmetric in $m,n$ and $k$. Therefore, there
exists a function $\F(\boldsymbol{t})$ such that
$$\frac{\pa^2 \F(\boldsymbol{t})}{\pa t_m\pa t_n} =
-mnb_{m,n}(\boldsymbol{t}).$$ By Proposition \ref{tau1}, the
conclusion  follows.
\end{proof}

\subsection{One variable reduction of dToda
hierarchy}\label{reduc}

\begin{proposition}\label{Pro2}Suppose that $(\mL(w;\boldsymbol{t}),
\tilde{\mL}(w;\boldsymbol{t}))$
is a solution of the dToda hierarchy whose dependence on
$\boldsymbol{t}=(\ldots, t_{-2}, t_{-1}, t_0, t_1, t_2, \ldots)$
is through a single variable $\lambda$. Namely, there exists a
function $\lambda(\boldsymbol{t})$ of $\boldsymbol{t}$ such that
\begin{align*}\mL(w;\boldsymbol{t})=&\mL(w,\lambda(\boldsymbol{t}))=
r(\lambda(\boldsymbol{t}))w+\sum_{n=0}^{\infty}
u_{n+1}(\lambda(\boldsymbol{t}))w^{-n}\\
\tilde{\mL}(w;\boldsymbol{t})^{-1}=
&\tilde{\mL}(w,\lambda(\boldsymbol{t}))^{-1}=r
(\lambda(\boldsymbol{t}))w^{-1}+\sum_{n=0}^{\infty}
\tilde{u}_{n+1}(\lambda(\boldsymbol{t}))w^{n}.\end{align*}
Then $\mL(w,\lambda)$ and
$\tilde{\mL}(w,\lambda)$
 satisfy the radial L\"owner
equations \eqref{radialout} and \eqref{radialin} respectively, in
which  $\phi(\lambda)=\log r(\lambda)$ and
$\eta(\lambda)=\sigma(\lambda)$. Moreover, if
$p(z,\lambda),\tilde{p}(z,\lambda)$ are the inverse functions of
$\mL(w,\lambda),\tilde{\mL}(w,\lambda)$ respectively;
$\Phi_n(w,\lambda), \Psi_n(w,\lambda), \, n\geq 1$ are Faber
polynomials of $p(z,\lambda)$ and $\tilde{p}(z,\lambda)$
respectively; $\xi_n(\lambda)$, $n\in\Z$ are defined by
\begin{align*}
\xi_n(\lambda)=\begin{cases}
\sigma(\lambda)\Phi_n'(\sigma(\lambda),\lambda),\hspace{1cm}&\text{if}\;\;\;
n\geq 1\\
1, &\text{if}\;\;\;n=0\\
\sigma(\lambda)\Psi_n'(\sigma(\lambda),\lambda),\hspace{1cm}&\text{if}\;\;\;
n\leq -1\end{cases};
\end{align*}
then $\lambda(\boldsymbol{t})$ satisfies the hydrodynamic type
equation
\begin{align}\label{hodograph2}
\frac{\pa\lambda}{\pa t_n}=\xi_n(\lambda)\frac{\pa\lambda}{\pa
t_0}.
\end{align}
\end{proposition}
\begin{proof}
Let $\F(\boldsymbol{t})=\log\tau_{\dToda}(\boldsymbol{t})$ be the
free energy of the solution. From the fourth equation in
\eqref{equation1}, we have
\begin{align*}
\left(\frac{1}{2}\pa_{t_0}+D(z_1)\right)\log p(z_2,\lambda)
=\left(\frac{1}{2}\pa_{t_0}+ D(z_2)\right)\log p(z_1,\lambda),
\end{align*}
Setting $z_1=z$ and tending $z_2 \to \infty$, we have
\begin{align*}
D(z)\phi= -\frac{1}{2}\pa_{t_0}\log \left(rp(z,\lambda)\right),
\hspace{1cm}\text{where}\;\;\;\phi= \log
r=\frac{1}{2}\frac{\pa^2\F}{\pa t_0^2}.
\end{align*}Suppose $\pa_{\lambda}\phi\neq 0$, then we have by the
chain rule:
\begin{align}\label{lambda1}
D(z)\lambda =- \frac{1}{2}\pa_{\lambda}\log
\left(rp(z,\lambda)\right)\frac{\pa_{t_0}\lambda}{\pa_{\lambda}\phi}
\end{align}
Applying $\pa_{t_0}$ to both sides of the first equation in
\eqref{equation1} and using the chain rule, we get
\begin{align*}
&\frac{p(z_1,\lambda)\pa_{\lambda}\log\left(rp(z_1,\lambda)\right)
-p(z_2,\lambda)\pa_{\lambda}\log\left(rp(z_2,\lambda)\right)}{p(z_1,\lambda)
-p(z_2,\lambda)}\pa_{t_0}\lambda\\
=&D(z_1)D(z_2)\pa_{t_0}\F=-D(z_1)\log\left(rp(z_2,\lambda)\right)\\
=&\frac{1}{2}\pa_{\lambda}
\log\left(rp(z_2,\lambda)\right)\pa_{\lambda}\log
\left(rp(z_1,\lambda)\right)\frac{\pa_{t_0}\lambda}{\pa_{\lambda}\phi}.
\end{align*}
If $\pa_{t_0}\lambda\neq 0$, then we can rewrite this in the form
\begin{align*}
&\Bigl(p(z_1,\lambda)\pa_{\lambda}\log\left(rp(z_1,\lambda)\right)
-p(z_2,\lambda)\pa_{\lambda}\log\left(rp(z_2,\lambda)\right)\Bigr)
\pa_{\lambda}\phi\\=&\frac{1}{2}\pa_{\lambda}
\log\left(rp(z_2,\lambda)\right)\pa_{\lambda}\log
\left(rp(z_1,\lambda)\right)(p(z_1,\lambda)-p(z_2,\lambda))
\end{align*}
which implies that
\begin{align}\label{definek}
\frac{1}{\sigma(\lambda)}=
\frac{1}{p(z,\lambda)}-
\frac{2\pa_{\lambda}\phi(\lambda)}{p(z,\lambda)\pa_{\lambda}\left(\phi+\log
p(z,\lambda)\right)}=-\frac{\pa_{\lambda}\phi-\pa_{\lambda}\log
p(z,\lambda)}{p(z,\lambda)\left(\pa_{\lambda}\phi+\pa_{\lambda}\log
p(z,\lambda)\right)}
\end{align}
is a constant independent of $z$. Rearranging, we obtain:
$$\frac{\pa p(z,\lambda)}{\pa\lambda}=
p(z,\lambda)
\frac{\sigma(\lambda)+p(z,\lambda)}{\sigma(\lambda)-p(z,\lambda)}
\frac{\pa\phi(\lambda)}{\pa\lambda},$$
which is the radial L\"owner equation \eqref{radialout} for $p(z,\lambda)$.

Similarly, from the last equation of \eqref{equation1}, we have
\begin{align}\label{lambda2}\tilde{D}(z)\lambda = \frac{1}{2}\pa_{\lambda}\log
\left(r\tilde{p}(z,\lambda)^{-1}\right)
\frac{\pa_{t_0}\lambda}{\pa_{\lambda}\phi}.\end{align}
Using
this and applying $\pa_{t_0}$ to the third equation of
\eqref{equation1}, we find that
\begin{align*}
\eta(\lambda)=\tilde{p}(z,\lambda)-\frac{2\tilde{p}(z,\lambda)
\pa_{\lambda}\phi}{\pa_{\lambda}\phi-\pa_{\lambda}\log
\tilde{p}(z,\lambda)}=-\tilde{p}(z,\lambda)\frac{\pa_{\lambda}
\phi+\pa_{\lambda}\log\tilde{p}(z,\lambda)}
{\pa_{\lambda}\phi-\pa_{\lambda}\log\tilde{p}(z,\lambda)}
\end{align*}
is a constant independent of $z$. This gives us
$$\frac{\pa \tilde{p}(z,\lambda)}{\pa\lambda}=
\tilde{p}(z,\lambda)\frac{\eta(\lambda)+\tilde{p}(z,\lambda)}
{\eta(\lambda)-\tilde{p}(z,\lambda)}\frac{\pa\phi(\lambda)}{\pa\lambda},
$$which is the radial L\"owner equation \eqref{radialin} for
$\tilde{p}(z,\lambda)$.

Now, applying $\pa_{t_0}$ to both sides of the second equation in
\eqref{equation1}, we have
\begin{align*}
&\frac{-\pa_{\lambda}\tilde{p}(z_2,\lambda)+p(z_1,\lambda)^{-1}
\tilde{p}(z_2,\lambda)
\pa_{\lambda}p(z_1,\lambda)}{p(z_1,\lambda)-\tilde{p}(z_2,\lambda)}
\pa_{t_0}\lambda=-\tilde{D}(z_2)\log
(rp(z_1,\lambda))\\
=&-\frac{1}{2}\pa_{\lambda}\log (rp(z_1,\lambda))\pa_{\lambda}\log
\left(r\tilde{p}(z_2,\lambda)^{-1}\right)
\frac{\pa_{t_0}\lambda}{\pa_{\lambda}\phi}.
\end{align*}Substituting the expressions for
$\pa_{\lambda}p(z_1,\lambda)$ and
$\pa_{\lambda}\tilde{p}(z_2,\lambda)$ obtained above, we get
\begin{align*}
-\frac{\eta(\lambda)+\tilde{p}(z_2,\lambda)}{\eta(\lambda)-
\tilde{p}(z_2,\lambda)}
+\frac{\sigma(\lambda)+p(z_1,\lambda)}{\sigma(\lambda)-p(z_1,\lambda)}=
\frac{2\sigma(\lambda)(p(z_1,\lambda)
-\tilde{p}(z_2,\lambda))}{
(\eta(\lambda)-\tilde{p}(z_2,\lambda))(\sigma(\lambda)-p(z_1,\lambda))},
\end{align*}
which after simplification implies that
$\sigma(\lambda)=\eta(\lambda)$.

Finally using the L\"owner equations for $p(z,\lambda)$ and
$\tilde{p}(z,\lambda)$ and applying \eqref{Fa1}, equations
\eqref{lambda1} and \eqref{lambda2} give
\begin{align*}
\sum_{n=1}^{\infty}\frac{z^{-n}}{n}\frac{\pa\lambda}{\pa t_n}=&
-\frac{\sigma(\lambda)}{\sigma(\lambda)-p(z_1,\lambda)}\frac{\pa\lambda}{\pa
t_0}=\sum_{n=1}^{\infty}\frac{\xi_n(\lambda)}{n}z^{-n}\frac{\pa\lambda}{\pa
t_0}\\
\sum_{n=1}^{\infty}\frac{z^n}{n}\frac{\pa\lambda}{\pa
t_{-n}}=&-\frac{\tilde{p}(z_2,\lambda)}{\sigma(\lambda)-\tilde{p}(z_2,
\lambda)}\frac{\pa\lambda}{\pa
t_0}
=\sum_{n=1}^{\infty}\frac{\xi_{-n}(\lambda)}{n}z^{n}\frac{\pa\lambda}{\pa
t_0},
\end{align*}
which imply
\begin{align*} \frac{\pa\lambda}{\pa
t_n}=\xi_n(\lambda)\frac{\pa\lambda}{\pa t_0}
\end{align*}
for all $n\in\Z$.
\end{proof}
\begin{remark}
If we assume that $\mL(w, \lambda(\boldsymbol{t}))$ maps the
exterior disc $\Del^*$ to a slit subdomain
$\sf{B}\setminus$$\Gamma_{\lambda}$ of a domain $\sf{B}$, then
since $\mL(\sigma(\lambda),\lambda)$ is  the tip of
$\Gamma_{\lambda}$, we have $|\sigma(\lambda)|=1$.
\end{remark}
\begin{remark}
Since $\eta(\lambda)=\sigma(\lambda)$, it appears that
$\mL(w,\lambda)$ and $\tilde{\mL}(w,\lambda)$ obey the same
differential equation. However, one should notice that
$\mL(w,\lambda)$ is a power series defined in a neighbourhood of
infinity and $\tilde{\mL}(w,\lambda)$ is a power series  defined
in a neighbourhood of the origin.
\end{remark}

Similarly to the dKP case, the hydrodynamic type equation
\eqref{hodograph2} can be solved by the general hodograph method
of Tsarev \cite{Tsarev}:
\begin{lemma}\label{le2}Let $R(\lambda)$ be any function of
$\lambda$. If $\lambda(\boldsymbol{t})$ is defined implicitly by
the hodograph relation
\begin{align}\label{ho3}t_0 +\sum_{n=1}^{\infty}\xi_n(\lambda) t_n
+\sum_{n=1}^{\infty}\xi_{-n}(\lambda)t_{-n}=R(\lambda),\end{align}
then $\lambda(\boldsymbol{t})$ satisfies \eqref{hodograph2}.
\end{lemma} \begin{remark}If for some $\lambda_0$,
$$\tilde{\mL}(w,\lambda_0)=\ov{\mL(1/\bar{w},\lambda_0)}^{\;-1},$$then
by uniqueness of solutions to differential equations,
\begin{align*}
\tilde{\mL}(w,\lambda)=\ov{\mL(1/\bar{w},\lambda)}^{\;-1}
\hspace{1cm}\text{for all}\;\lambda.
\end{align*}In this case,
$$w\Psi'_n(w,\lambda)=-\ov{\frac{1}{\bar{w}}\Phi_n'\left(\frac{1}{\bar{w}},\lambda\right)}.$$If
we also impose the condition $|\sigma(\lambda)|=1$, then
$$\xi_{-n}(\lambda)=-\ov{\xi_n(\lambda)}.$$ Therefore, if $t_0$
is a real-valued variable
and  for $n\neq 0$, $t_n$ are complex-valued variables such that
$t_{-n}=-\bar{t}_n$, then \eqref{ho3} defines
$\lambda(\boldsymbol{t})$ as a real-valued function in a certain
domain of $\boldsymbol{t}$.
\end{remark}
As a converse to Proposition \ref{Pro2}, we have
\begin{proposition}\label{Pr2}
Suppose
$$G(w,\lambda)=e^{\phi(\lambda)}w+\sum_{n=0}^{\infty}u_{n+1}(\lambda) w^{-n}$$
is a
solution of the radial L\"owner equation \eqref{radialout} and
\begin{align*} F(w,\lambda)=\ov{G(1/\bar{w},\lambda)}^{\;-1}.
\end{align*}Let
$\lambda(\boldsymbol{t})$ be a solution of equation
\eqref{hodograph2}, where $\boldsymbol{t}=\{t_n\}_{n\in\Z}$, $t_0$
is a real-valued variable and for $n\neq 0$,  $t_n$ are complex
variables such that $t_{-n}=-\bar{t}_n$. Under these conditions,
$(\mL(w;\boldsymbol{t}), \tilde{\mL}(w; \boldsymbol{t}))$ defined
by
$$\mL(w;\boldsymbol{t})=G(w,\lambda(\boldsymbol{t}))\hspace{1cm}
\text{and}\hspace{1cm}
\tilde{\mL}(w;\boldsymbol{t})=F(w,\lambda(\boldsymbol{t}))$$is a
solution of the dToda hierarchy \eqref{Lax1} with
$t_{-n}=-\bar{t}_n$, $n\geq 1$.
\end{proposition}
\begin{proof}
First, it is easy to verify that $F(w,\lambda)$ satisfies the
radial L\"owner equation \eqref{radialin} with
$\eta(\lambda)=\sigma(\lambda)$.

Let $g(z,\lambda)$ and $f(z,\lambda)$ be the inverse functions of
$G(w,\lambda)$ and $F(w,\lambda)$ respectively and let
$b_{m,n}(\boldsymbol{t})=b_{m,n}(\lambda(\boldsymbol{t}))$,
$m,n\in \Z$ be the Grunsky coefficients of $(f(z,\lambda),
g(z,\lambda))$. Applying $\pa_{t_k}$ to both sides of the first
equation in \eqref{Grunsky} and using the L\"owner equation for
$g(z,\lambda)$, equation \eqref{hodograph2} and the first equation
in \eqref{Fa1}, we have
\begin{align*}
-\sum_{m=1}^{\infty} \sum_{n=1}^{\infty}\frac{\pa
b_{m,n}(\boldsymbol{t})}{\pa t_k}z_1^{-m}z_2^{-n} =&
\frac{\pa_{\lambda}\left(e^{\phi(\lambda)}
g(z_1,\lambda)\right)-\pa_{\lambda}\left(e^{\phi(\lambda)}g(z_2,\lambda)\right)}
{e^{\phi(\lambda)}\left(g(z_1,\lambda)-g(z_2,\lambda)\right)}\pa_{t_k}\lambda\\
=&\frac{2\sigma(\lambda)^2\pa_{\lambda}\phi\xi_k(\lambda)\pa_{t_0}\lambda}
{(\sigma(\lambda)-g(z_1,\lambda))(\sigma(\lambda)-g(z_2,\lambda))}\\
=&2\sum_{m=1}^{\infty}\sum_{n=1}^{\infty}\frac{1}{mn}\xi_m(\lambda)\xi_{n}(\lambda)z_1^{-m}z_2^{-n}
\xi_k(\lambda)\pa_{\lambda}\phi\pa_{t_0}\lambda.
\end{align*}Similarly, applying $\pa_{t_k}$ to both sides of
equations in \eqref{Grunsky2} and \eqref{Grunsky3} and using the
L\"owner equations for $g(z,\lambda)$ and $f(z,\lambda)$, Lemma
\ref{le2} and equations \eqref{Fa1}, we obtain
\begin{align*}
-\sum_{m=1}^{\infty} \sum_{n=1}^{\infty}\frac{\pa
b_{m,-n}(\boldsymbol{t})}{\pa
t_k}z_1^{-m}z_2^{n}=&2\sum_{m=1}^{\infty}\sum_{n=1}^{\infty}\frac{1}{mn}
\xi_m(\lambda)\xi_{-n}(\lambda)z_1^{-m}z_2^{n}
\xi_k(\lambda)\pa_{\lambda}\phi\pa_{t_0}\lambda\\
-\sum_{m=1}^{\infty} \sum_{n=1}^{\infty}\frac{\pa
b_{-m,-n}(\boldsymbol{t})}{\pa
t_k}z_1^{m}z_2^{n}=&2\sum_{m=1}^{\infty}\sum_{n=1}^{\infty}\frac{1}{mn}
\xi_{-m}(\lambda)\xi_{-n}(\lambda)z_1^{m}z_2^{n}
\xi_k(\lambda)\pa_{\lambda}\phi\pa_{t_0}\lambda\\
-\sum_{n=1}^{\infty}\frac{\pa b_{n,0}(\boldsymbol{t})}{\pa t_k}
z^{-n}=&-2\sum_{n=1}^{\infty}\frac{\xi_n(\lambda)}{n}z^{-n}\xi_k(\lambda)\pa_{\lambda}\phi\pa_{t_0}\lambda\\
-\sum_{n=1}^{\infty}\frac{\pa b_{-n,0}(\boldsymbol{t})}{\pa t_k}
z^{n}=&-2\sum_{n=1}^{\infty}\frac{\xi_{-n}(\lambda)}{n}z^{n}\xi_k(\lambda)\pa_{\lambda}\phi\pa_{t_0}\lambda.
\end{align*}Comparing coefficients, we find that
\begin{align*}
2\xi_m(\lambda)\xi_{n}(\lambda)\xi_k(\lambda)\pa_{\lambda}\phi\pa_{t_0}\lambda
=\begin{cases} -|mn|\pa_{t_k} b_{m,n}(\boldsymbol{t}),
\hspace{1cm}&\text{if}\;\;|m|\neq 0, n\neq 0\\
|m|\pa_{t_k} b_{m,0}(\boldsymbol{t}), &\text{if}\;\;m\neq 0 ,n= 0\\
2\pa_{t_k}\phi(\lambda(\boldsymbol{t})),&\text{if}\;\;m=n=0.
\end{cases}
\end{align*}
Since the left hand side is completely symmetric in $m,n,k$,
we conclude that there exists a function $\F(\boldsymbol{t})$ such
that
\begin{align*}
\frac{\pa^2\F(\boldsymbol{t})}{\pa t_m\pa t_n}
=\begin{cases}-|mn|b_{m,n}(\boldsymbol{t}),
\hspace{1cm}&\text{if}\;\;m\neq 0, n\neq 0\\
|m| b_{m,0}(\boldsymbol{t}), &\text{if}\;\;m\neq 0,n= 0\\
2\phi(\lambda(\boldsymbol{t})),&\text{if}\;\;m=n=0.
\end{cases}
\end{align*}
Since $\phi(\lambda)=-b_{0,0}(\lambda)$, the assertion
follows from Proposition \ref{tau2}.
\end{proof}
\begin{remark}
In \cite{Manas}, Manas has considered the reduction of $r$-th
dispersionless Dym (dDym) hierarchy. When $r=1$, the dDym
hierarchy is gauge equivalent to 'half' of the dToda hierarchy
\eqref{Lax1}, where only the power series $\mL(w,\boldsymbol{t})$
depending on $\boldsymbol{t}=(t_0, t_1, t_2,\ldots )$\footnote{The
variable $t_0$ should be identified with $x$ in \cite{Manas}.} is
present. Ma{\~n}as proved a result similar to Proposition \ref{Pr2} by
a completely different method. 

Recently, Takasaki and Takebe \cite{TTnew} also showed that a solution
to the radial L\"owner equation together with the hydrodynamic type
equation \eqref{hodograph2} for $n\geq 1$ generates a solution to the
first series of the Lax equations of the dToda hierarchy (the first
equation in \eqref{Lax1}). In fact, the results of \cite{TTnew} are
easily recovered from \cite{Manas} by just adding linear and constant
terms in $p$ to reduction condition (30) of \cite{Manas}.
\end{remark}

\section{Large $N$ eigenvalue integrals and reductions of
dispersionless hierarchies}

Here we explain the relations between dispersionless integrable
hierarchies, their reductions and conformal maps from the
perspective of the model of normal random matrices
\cite{CY,ref20}. Actually, we need only the eigenvalue integral
for that model, which represents the partition function of the 2D
Coulomb gas in external field. In this section, the exposition is
on the physical level of rigor.

\subsection{Large $N$ integrals and dToda hierarchy}

A convenient starting point is the $N$-fold integral
\begin{equation}\label{tauN}
\tau_N = \frac{1}{N!} \int_{\CC^N}
\prod_{i<j}^{N} |z_i - z_j |^2
\prod_{k=1}^{N} e^{\frac{1}{\hbar} \sum_{n= 1}^{\infty}
(t_n z_{k}^{n}+\bar t_n \bar z_{k}^{n})} \, d\mu (z_k , \bar z_k)
\end{equation}
where $d\mu$ is some integration measure in the
complex plane, and $\hbar$ is a parameter. For $d\mu =
e^{\frac{1}{\hbar}W(z, \bar z)} d^2z$, $\tau_N$ is the partition
function of the model of normal random matrices with the potential
$W$ written as an integral over eigenvalues.
Clearly, $\tau_N$ (\ref{tauN})
has the meaning of the partition function for a system of $N$
2D Coulomb charges in an external potential.
 The basic fact about
the integral (\ref{tauN}) is:
\begin{itemize}
\item For any measure $d\mu$ (including singular measures
supported on sets of dimension less than $2$), $\tau_N$, as a
function of $t_0 =N\hbar$,
$\{t_n \}$,  $\{-\bar t_n \}$, is a $\tau$-function of
the (dispersionful) 2D Toda hierarchy with $t_{-n}=-\bar{t}_n$.
\end{itemize}
In a slightly different form, this statement first appeared in
\cite{ref19}, see also \cite{ref20}.

In the large $N$ limit ($N\to \infty$, $\hbar \to 0$, $t_0 = \hbar
N$ finite) $\tau_N$ generates the ``free energy'' $\mathcal{F}$
for the dToda hierarchy via
\begin{equation*}\label{Ftau}
\mathcal{F}(\boldsymbol{t})= \lim_{N\to \infty}
\left ( \hbar^2 \log \tau_N \right). \end{equation*}
where $\boldsymbol{t}= \{ \ldots , -\bar t_2 , -\bar t_1 ,
t_0 , t_1 , t_2 , \ldots \}$.
It obeys the
dispersionless Hirota equations \eqref{equation1}. Second order
derivatives of $\mathcal{ F}$ enjoy a nice geometric
interpretation through conformal maps. In short, this goes as follows.
As $N \to \infty$, the integral in (\ref{tauN})
is determined by the most
favorable configuration of $z_i$'s, i.e., the one at which the
integrand has a maximum (which becomes very sharp
in the large $N$ limit). In other words, the free energy
$\mathcal{F}$ is essentially the electrostatic energy of the
equilibrium configuration of the 2D Coulomb charges in the
external potential.
Exploiting further the electrostatic analogy, one
can see that in the equilibrium the ``charges''
densely fill a bounded domain $\DD$ in the complex plane. For
simplicity, we assume that $\DD$ is connected. In the matrix model
interpretation, this domain is called the support of eigenvalues.
Clearly, it is a compact subset of the support of the measure $d\mu$.

Let $p(z)$ be the conformal mapping function from the {\it exterior}
of the domain $\DD$ onto the exterior of the unit circle
normalized as $p(z)=z/r +O(1)$ at large $|z|$ with a real $r$
called the exterior conformal radius of the domain $\DD$. In
\cite{ref6_1,ref6_2,ref8,ref9} it was shown that the function
$p(z)$ can be expressed through ${\mathcal{F}}$ in the following
different but equivalent ways:
\begin{align}\label{confmap1}
rp(z)= \, &z\, e^{-D(z)\pa_{t_0}{\mathcal{F}}},\\ rp(z)= \, &z- a
-D(z)\pa_{t_1}{\mathcal{F}},\nonumber\\  rp^{-1}(z)= \, & D(z)
\pa_{\bar t_1}{\mathcal{F}} \nonumber\end{align} where
\begin{equation*}\label{r1} 2\log r =\frac{\pa^2  {\mathcal{F}}}{\pa
t_{0}^{2}}\,, \quad a =\frac{\pa^2 {\mathcal{F}}}{\pa t_{0} \pa
t_1}
\end{equation*} and the operator $D(z)$ is defined in \eqref{op}.
Moreover, the free energy obeys the Hirota equations
\eqref{equation1} with $\tilde p(z)=1/\overline{p(1/\bar z)}$.
It is convenient to rewrite them in terms of the function
$$\bar{p}(z) = \ov{p(\z)}= r^{-1}z\, e^{-\bar
D(z)\pa_{t_0}{\mathcal{F}}}\hspace{1cm}\text{where}\hspace{1cm}
\bar D(z)= \sum_{n=1}^{\infty} \frac{z^{-n}}{n}\,
\frac{\pa}{\pa \bar t_n}
$$
We have:
\begin{align}\label{Toda1} D(z_1)D(z_2){\mathcal{F}}=&\log
\frac{rp(z_1)-rp(z_2)}{z_1 -z_2} ,\\\bar D(z_1)\bar
D(z_2){\mathcal{F}}=&\log \frac{r\bar p(z_1)-r\bar p(z_2)}{z_1
-z_2}\nonumber\\  -  D(z_1)\bar D(z_2)\mathcal{ F}=&\log \left (
1-\frac{1}{ p(z_1)\bar p(z_2)}\right )\nonumber
\end{align}
We emphasize that all these relations hold
true for {\it any} measure $d\mu$ in (\ref{tauN}) provided the
equilibrium configuration of $z_i$'s at $N \to \infty$ is well
defined.

\begin{figure}[h]
 \begin{center}
  \psfrag{D}{$\DD$}
  \psfrag{Ga=supp(dmu)}{$\Gamma=\text{supp}(d\mu)$}
  \includegraphics{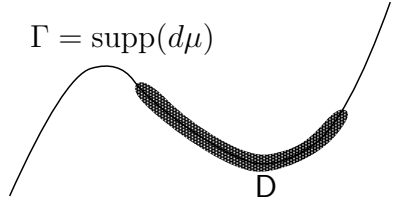}
 \end{center}
\label{fig:6-1-1}
\caption{The supports of the measure and of eigenvalues.}
\end{figure}

In particular, one may consider the measure supported on a curve
$\Gamma$, then the integral (\ref{tauN}) becomes one-dimensional
(along $\Gamma$) in each variable: \begin{equation}\label{tauN1}
\tau_N = \frac{1}{N!} \int_{\Gamma^N}
\prod_{i<j}^{N} |z_i - z_j |^2
\prod_{k=1}^{N} e^{\frac{1}{\hbar} \sum_{n\geq 1} (t_n
z_{k}^{n}+\bar t_n \bar z_{k}^{n})} \, |dz_k|
\end{equation}
In the large $N$ limit, the support of eigenvalues, $\DD$, is then an
arc of the curve $\Gamma$ (or several disconnected arcs, but we do
not consider this case in the present paper). The function $p(z)$
maps the slit domain $\C\setminus\complex{D}$ onto the exterior of
the unit circle. (See Figure~\ref{fig:6-1-1}.)

For the dToda hierarchy, the choice of the measure supported on a
curve means {\it a reduction}. A familiar example is the dToda
chain, where one may take $\Gamma$ to be either real or imaginary
axis. Consider a general (continuous) curve $\Gamma$ infinite in
both directions. It is clear that $p(z)$ and the Lax functions
(the functions inverse to $p(z)$ and $1/\bar p (1/z)$) depend on the
times through two parameters only. One can set them to be, for
example, the positions of the two ends of the arc $\DD$ on the
curve $\Gamma$ (although this is not a best choice). Reductions
with only one independent parameter are also possible. The
simplest possibility to obtain them from the large $N$ integral
(\ref{tauN}) is to take the measure $d\mu$ supported on a
half-infinite curve, starting at $0$ (for example) and going to
infinity, and such that the arc $\DD$ starts from $0$ as the times
independently vary in some open set. A more general class of
examples can be obtained in a similar way if one restricts the
measure to be supported on the boundary of some domain
${\sf B}$ and on a curve $\Gamma$ starting on the boundary
and coming to infinity. Under certain conditions,
only the second edge of the arc moves, so that we are left with
one real parameter $\lambda = \lambda (\boldsymbol{t})$.
We have shown in Section \ref{reduc} that in this case, the
corresponding conformal map $p(z)$ satisfies the radial L\"owner
equation.

\subsection{The reflection symmetry case:
the chordal L\"owner equation from reductions of the dToda
hierarchy}

The case when the measure $d\mu$ has a symmetry
is of special interest. In order to respect the symmetry, one may need to
restrict the flows of the full dToda hierarchy to a submanifold of
the space  with coordinates $t_n$, $\bar t_n$. Reductions of the so obtained
sub-hierarchy may lead to interesting classes of conformal maps.
Below we study the case of reflection symmetry.

Suppose the measure is symmetric under the reflection in the real
axis, i.e., $d\mu (z, \bar z)= d\mu (\bar z , z)$. We will show
that this case allows one to relate conformal maps described
by the chordal L\"owner equation \eqref{chordal} to the conformal
maps generated by solutions of the dToda hierarchy. The relation
is not immediately obvious.

The symmetry of the measure does not yet imply the symmetry of the
support of eigenvalues since the latter depends also on the
function $\sum_n (t_n z^n +\bar t_n \bar z^n)$ which does not
enjoy the reflection symmetry unless all $t_k$'s are real.
If they are, then the domain $\DD$ is invariant under
the reflection (complex conjugation). In this case the conformal
map has real coefficients: $\bar p(z)=p(z)$, i.e.,
$\pa_{t_0}\pa _{t_n}\mathcal{F}=\pa_{t_0}\pa _{\bar t_n}\mathcal{F}$.
Equations \eqref{Toda1}
then imply that in the real section of the times manifold ($t_n
=\bar t_n$) it holds
$\pa_{t_m}\pa _{t_n}\mathcal{F}=\pa_{\bar t_m}\pa _{\bar t_n}\mathcal{F}$
and
$\pa_{t_m}\pa _{\bar t_n}\mathcal{F}=\pa_{\bar t_m}\pa _{t_n}\mathcal{F}$
for any $m,n$.
Therefore, the dToda hierarchy restricted to this real section
describes conformal maps of symmetric domains. Unfortunately, we
do not know whether the restricted flows form a reasonable
hierarchy.

Deformations which destroy the reflection symmetry of $\DD$ can be
parametrized by the times $s_n =\frac{1}{\sqrt{2}}(t_n - \bar
t_n)$ (we define them to be purely imaginary for later
convenience). We are going to show that infinitesimal deformations
of this kind are described by the dKP hierarchy in the times
$\boldsymbol{s}=\{s_n \}$. More precisely, the Hirota equation of the form
\eqref{equation2} generating the dKP hierarchy in the times
$\boldsymbol{s}$
holds true at the point $\boldsymbol{s} =0$,
and for $s_n$ of order $\epsilon$ the
corrections to this equation are of order $\epsilon ^2$.

\begin{figure}[h]
 \begin{center}
  \psfrag{D}{$\DD$}
  \psfrag{Ga=supp(dmu)}{$\Gamma=\text{supp}(d\mu)$}
  \includegraphics{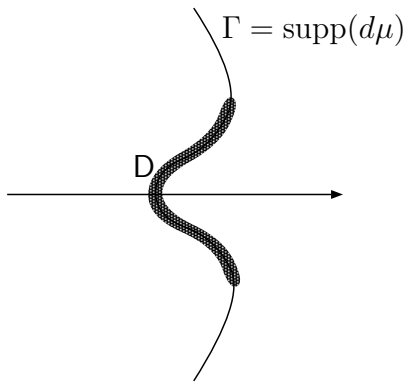}
 \end{center}
\label{fig:6-2-1}
\caption{Reflection symmetry case; a slit domain.}
\end{figure}

Let us take the sum of the three equations \eqref{Toda1} and the
fourth equation obtained from the third one in \eqref{Toda1}
by the interchange $z_1 \leftrightarrow z_2$. After
exponentiating, we get, using the symmetry $\bar p(z)=p(z)$:
$$
(z_1 \!-\! z_2)e^{\frac{1}{2} (D(z_1)\! -\! \bar D(z_1)) (D(z_2)\!
-\! \bar D(z_2))\mathcal{F}} =r\left ( p(z_1) \! - \! p(z_2) \!
+\! p^{-1}(z_1) \! - \! p^{-1}(z_2)\right )
$$
Summing three equations of this type written for all pairs of the
points $z_1 , z_2 , z_3$, we obtain:
\begin{equation}\label{dKPs} (z_1 \!-\! z_2)e^{D^{\rm
KP}(z_1)D^{\rm KP}(z_2)\mathcal{F}} +(z_2 \!-\! z_3)e^{D^{\rm
KP}(z_2)D^{\rm KP}(z_3)\mathcal{F}} +(z_3 \!-\! z_1)e^{D^{\rm
KP}(z_3)D^{\rm KP}(z_1)\mathcal{F}}=0 \end{equation} where
\begin{equation*}\label{refl1} D^{\rm KP}(z)=
\sum_{n=1}^{\infty}\frac{z^{-n}}{n}
\pa_{s_n}\,, \quad
\pa_{s_n}=\frac{1}{\sqrt{2}}(\pa_{t_n}-\pa_{\bar t_n})
\end{equation*}Setting
\begin{equation}\label{refl22} k(z)=z-D^{{\rm
KP}}(z)\pa_{s_1}\mathcal{F} \end{equation} and letting
$z_3\rightarrow \infty$ in \eqref{dKPs}, we get
\begin{equation}\label{dKPs1} D^{\rm KP}(z_1)D^{\rm
KP}(z_2)\mathcal{F} =\log \frac{k(z_1)-k(z_2)}{z_1 - z_2},
\end{equation} which is the Hirota equation for the dKP hierarchy
\eqref{equation2}. 
It is important to note that equations
(\ref{dKPs}), (\ref{dKPs1}) {\it hold at the point
$\boldsymbol{s}=0$ only}
and thus the evolution in $s_n$ that they define is not consistent
with the one defined by the dToda hierarchy for the same function
$\mathcal{F}$.
However, for small $s_n$ the two evolutions agree up to
higher order terms. More precisely, let all the $s_n$'s be
of order $\epsilon \to 0$, then the coefficients
of $p(z)$ acquire imaginary parts of the same order $\epsilon$
(or higher) while their real parts are not changed up to this order.
Summing the equations \eqref{Toda1} as before but
without the assumption that $\bar p(z)=p(z)$, one can see that
the terms of order $\epsilon$ in the r.h.s. cancel.
Therefore,
the l.h.s. of (\ref{dKPs}) is of order $\epsilon^2$.
This means that the
equation obtained from (\ref{dKPs}) by applying any derivative
$\pa_{s_n}$ at $s_n =0$ is still valid. The same is true for
equation (\ref{dKPs1}). As we have seen, this is enough to derive
the chordal L\"owner equation from the dKP hierarchy (with a
reduction imposed). Therefore, conformal maps of slit domains {\it
symmetric under reflection in the real axis} satisfy the chordal
L\"owner equation \eqref{chordal}, as they should. Note that
symmetric slits have only one parameter (the curve $\Gamma$ is
fixed). (See Figure~\ref{fig:6-2-1}.)

It remains to verify that the function $k(z)$ defined by
(\ref{refl22}) is indeed the conformal map with the required
properties\footnote{In the L\"owner equation literature, the map
$k(z)$ is usually regarded as a map from the upper half plane with
a slit onto the upper half plane. The Schwarz symmetry principle
allows one to analytically continue the map to the lower half
plane, with the result being characterized by these properties.}:
\begin{itemize}
\item[(i)] Normalization $k(z)= z + O(1/z)$ as $z \to \infty$;
\item[(ii)] Real coefficients;
\item[(iii)] The image is the
complement of a segment of the real axis.
\end{itemize}
The first two are obvious from the representation (\ref{refl22}).
The third one follows from the identity
\begin{equation}\label{refl3}
k(z)=r\left (p(z) + p^{-1}(z) \right )
+a \end{equation}
which is a direct consequense of the second
and third equations in \eqref{confmap1}. Indeed, the function
$r(w+w^{-1})+a$ maps the exterior of the unit circle in the $w$-plane,
which is the image of the map $p(z)$,
 onto the exterior of a segment of the real line.
The function $k(z)$ is thus the composition of the two maps.

\vspace{1cm} \noindent \textbf{Acknowledgement} \;\; 
The work of T.T. was supported in part by JSPS grant 15540014 and by
JSPS/RFBR joint project. The work of L-P.T. was supported in part by
MMU funding PR/2006/0590. The work of A.Z. was supported in part by RFBR
grant 04-01-00642, by grant NSh-8004.2006.2 for support of scientific
schools and by grant NWO 047.017.015. We would like to thank
E. Bettelheim, I. Gruzberg, A. Marshakov, S. Natanzon, K. Takasaki,
P. Wiegmann for valuable discussions.

\appendix
\section{Examples}
\subsection{Chordal L\"owner equation and dKP hierarchy}

\subsubsection{}\textbf{Example 1}\\
We consider the Chordal L\"owner equation \eqref{chordal} with
$a_1=-\lambda$, $U(\lambda)=U\in\R$ a constant and initial
condition $H(w,0)=w$. It is easy to find that its solution is
given by
\begin{align}\label{so1}
H(w,\lambda)=&U+w\sqrt{1-2Uw^{-1}+(U^2-2\lambda)w^{-2}}=
w-\frac{\lambda}{w}-\frac{\lambda U}{w^2}+\ldots\\
h(z,\lambda)=&U+z\sqrt{1-2Uz^{-1}+(U^2+2\lambda)z^{-2}}=
z+\frac{\lambda}{z}+\frac{\lambda
U}{z^2}+\ldots\nonumber
\end{align}$H(w,\lambda)$ maps the upper half plane $\U$ conformally onto
 $\U\setminus\Gamma_{\lambda}$, where $\Gamma_{\lambda}$
is the line segment $$\Gamma_{\lambda}= \{ U+i\alpha\;:\;
\alpha\in[0, \sqrt{2\lambda}]\}.$$ Let $\Phi_n(w,\lambda)$, $n\geq
1$, be the Faber polynomials of $h(z,\lambda)$. A straight\-for\-ward
computation gives
\begin{align*}
\Phi_1(w,\lambda)=w,
\hspace{1cm}\Phi_2(w,\lambda)=w^2-2\lambda,\hspace{1cm}\Phi_3(w,\lambda)=w^3-3\lambda
w-3\lambda U.
\end{align*}Therefore the functions $\chi_n(\lambda)$, $n=1,2,3$ in
\eqref{hodograph} are given by$$\chi_1(\lambda)=1,
\hspace{1cm}\chi_2(\lambda)=2U,
\hspace{1cm}\chi_3(\lambda)=3U^2-3\lambda.$$When $t_n=0$ for all
$n\geq 4$, the hodograph relation \eqref{ho2} with $R(\lambda)=0$
reads as
\begin{align*}
x+t_1+2Ut_2+(3U^2-3\lambda)t_3=0
\end{align*}
which for $\{x+t_1+2Ut_2+3U^2t_3\geq 0 $ and $t_3>0\}$ or
$\{x+t_1+2Ut_2+3U^2t_3\leq 0 $ and $t_3<0\}$ gives
\begin{align*} \lambda=\frac{x+t_1+2Ut_2+3U^2t_3}{3t_3}\geq 0.
\end{align*}Substituting into the first equation of \eqref{so1},
we find that the
power series
\begin{align*}
\left.\mL(w; \boldsymbol{t})\right|_{t_n=0 \, (n\geq
4)}=U+w\sqrt{1-2Uw^{-1}-\left(U^2+\frac{2(x+t_1+2Ut_2)}{3t_3}\right)w^{-2}}
\end{align*}is a solution of the dKP hierarchy \eqref{Lax2}.

\subsubsection{}\textbf{Example 2}\\
Let \begin{align}\label{soso}
H(w,\lambda)=w+\frac{\lambda^2}{w-2\lambda}=w+\lambda^2
w^{-1}+2\lambda^3w^{-2}+4\lambda^4 w^{-3}+\ldots.
\end{align}For $\lambda>0$, $H(w,\lambda)$ maps the upper half plane conformally
onto $$\C\setminus (\{x\; :\; x\in (-\infty,0]\}\cup\{x\;:\;x\in [
4\lambda,\infty)\}).$$It is easy to check that $H(w,\lambda)$
satisfies the chordal L\"owner equation \eqref{chordal} with
$a_1=\lambda^2$ and $U(\lambda)=3\lambda$. Let
$\Phi_n(w,\lambda)$, $n\geq 1$ be the Faber polynomials of
$h(z,\lambda)$, the inverse function of $H(w,\lambda)$, then
\begin{align*}
\Phi_1(w,\lambda)=w,
\hspace{1cm}\Phi_2(w,\lambda)=w^2+2\lambda^2,\hspace{1cm}
\Phi_3(w,\lambda)=w^3+3\lambda^2w+6\lambda^3.
\end{align*}The functions $\chi_n(\lambda)$, $n=1,2,3$ defined  in
\eqref{hodograph} are given by
\begin{align*}
\chi_1(\lambda)= 1,\hspace{1cm}\chi_2(\lambda)=6\lambda,
\hspace{1cm} \chi_3(\lambda)= 30\lambda^2.
\end{align*}When $t_n=0$ for all
$n\geq 3$, the hodograph relation \eqref{ho2} with $R(\lambda)=0$
reads as
\begin{align*}
x+t_1+6t_2\lambda=0,
\end{align*}
which for $\{x+t_1\geq 0$ and $t_2<0\}$ or $\{x+t_1\leq 0$ and
$t_2>0\}$ gives
\begin{align*} \lambda=-\frac{x+t_1}{6t_2}\geq 0.
\end{align*}Substituting into  \eqref{soso}, we find that the
power series
\begin{align*}
\left.\mL(w, \boldsymbol{t})\right|_{t_n=0 \, (n\geq
3)}=w+\frac{(x+t_1)^2}{12 t_2}\frac{1}{3t_2w+(x+t_1)}
\end{align*}is a solution of the dKP hierarchy \eqref{Lax2}.
When $t_n=0$ for all $n\geq 4$, the hodograph relation \eqref{ho2}
with $R(\lambda)=0$ reads as
\begin{align*}
x+t_1+6t_2\lambda+30t_3\lambda^2=0,
\end{align*}which can be solved for
$\left.\lambda(\boldsymbol{t})\right|_{t_n=0 \, (n\geq 4)}$ in a
certain domain of $\boldsymbol{t}$.

\subsection{Radial L\"owner equation and dToda hierarchy}

\subsubsection{}\textbf{Example 1}

 We consider the radial L\"owner equation
\eqref{radialout} with $\phi(\lambda)=\lambda$,
$\sigma(\lambda)=\sigma\in S^1$ a constant and initial condition
$G(w,0)=w$. It is easy to check that its solution is given by
\footnote{In fact, up to a re-parametrization of $\lambda$, this
example is the same as the one given in \cite{TTnew}.}
\begin{align}\label{so2}
G(w,\lambda)=&-\sigma+\frac{e^{\lambda}w}{2}\left(\left(1+\frac{\sigma}{w}\right)^2
+\left(1+\frac{\sigma}{w}\right)
\sqrt{1+\frac{2\sigma(1-2e^{-\lambda})}{w}+\frac{\sigma^2}{w^{2}}}\right)\\
=&e^{\lambda}w+2\sigma(e^{\lambda}-1)+\frac{\sigma^2(e^{\lambda}-e^{-\lambda})}{w}+
\frac{2\sigma^3e^{-\lambda}(1-e^{-\lambda})}{w^2}+\ldots\nonumber\\
g(z,\lambda)=&-\sigma+\frac{e^{-\lambda}z}{2}\left(\left(1+\frac{\sigma}{z}\right)^2+
\left(1+\frac{\sigma}{z}\right)\sqrt{1+\frac{2\sigma(1-2e^{\lambda})}{z}+\frac{\sigma^2}{z^{2}}}\right).\nonumber
\end{align}$G(w,\lambda)$ maps the exterior disc $\Del^*$
conformally onto $\Del^*\setminus\Gamma_{\lambda}$, where
$\Gamma_{\lambda}$ is the line segment
\begin{align*}
\Gamma_{\lambda}=\{ \alpha\sigma\;:\;\alpha\in[1,
2e^{\lambda}-1+2\sqrt{e^{2\lambda}-e^{\lambda}}].\}
\end{align*}The maps
\begin{align}\label{so3}
F(w,\lambda)
=&-\sigma+\frac{e^{\lambda}}{2w}\left((w+\sigma)^2-\sigma(w+\sigma)\sqrt{1+2\sigma^{-1}(1-2e^{-\lambda})w+\sigma^{-2}w^2}
\right)\\
f(z,\lambda)
=&-\sigma+\frac{e^{-\lambda}}{2z}\left((z+\sigma)^2-\sigma(z+\sigma)\sqrt{1+2\sigma^{-1}(1-2e^{\lambda})z+\sigma^{-2}z^2}
\right)\nonumber
\end{align}
satisfy the radial L\"owner equations \eqref{radialin} with
$\eta(\lambda)=\sigma$ and initial condition $F(w,0)=w$. It is
easy to verify directly that
\begin{align}\label{inver}
F(w,\lambda)=\ov{G(1/\bar{w},\lambda)}^{\;-1}.\end{align}
Therefore, $F(w,\lambda)$ maps the unit disc $\Del$ conformally
onto $\Del\setminus\tilde{\Gamma}_{\lambda}$, where
$\tilde{\Gamma}_{\lambda}$ is the line segment
\begin{align*}
\tilde{\Gamma}_{\lambda}=\{ \alpha\sigma\;:\;\alpha\in[
2e^{\lambda}-1-2\sqrt{e^{2\lambda}-e^{\lambda}},1].\}
\end{align*}Let $J_{\lambda}=g(\Gamma_{\lambda}, \lambda)$.
In fact, for fixed $\lambda$, $F(w,\lambda)$ is the analytic
continuation of $G(w,\lambda)$ to the set $\Del\cup (S^1\setminus
 J_{\lambda})$.

 Let $\Phi_n(w,\lambda)$, $\Psi_n(w,\lambda)$, $n\geq 1$ be the
Faber polynomials of $g_{\lambda}(z)$ and $f_{\lambda}(z)$
respectively. Then
\begin{align*}
\Phi_1(w,\lambda)=&e^{\lambda}w
+2\sigma(e^{\lambda}-1),\hspace{1cm}
\Psi_1(w,\lambda)=e^{\lambda}w^{-1}+2\sigma^{-1}(e^{\lambda}-1),\\
\Phi_2(w,\lambda)=&e^{2\lambda}w^2+4\sigma
e^{\lambda}(e^{\lambda}-1)w+2\sigma^2(e^{\lambda}-1)(2e^{\lambda}-1+e^{-\lambda}),\\
\Psi_2(w,\lambda)=&e^{2\lambda}w^{-2}+4\sigma^{-1}e^{\lambda}(e^{\lambda}-1)w^{-1}
+2\sigma^{-2}(e^{\lambda}-1)(2e^{\lambda}-1+e^{-\lambda}).
\end{align*}Therefore, the functions $\xi_n(\lambda)$, $n=\pm 1, \pm 2$ in \eqref{hodograph2} are given by
\begin{align*}
\xi_1(\lambda)=& \sigma e^{\lambda}, \hspace{3cm}
\xi_{-1}(\lambda) =-\sigma^{-1}e^{\lambda}\\
\xi_2(\lambda)=&2\sigma^2e^{\lambda}(3e^{\lambda}-2),
\hspace{1cm}\xi_{-2}(\lambda)=
-2\sigma^{-2}e^{\lambda}(3e^{\lambda}-2).
\end{align*}When $t_{-1}=-\bar{t}_1$ and $t_n=0$ for $|n|\geq 2$, the hodograph relation
\eqref{ho3} with $R(\lambda)=0$ reads as
\begin{align}\label{sol1}
t_0+t_1 \sigma e^{\lambda}+\bar{t}_1\bar{\sigma}e^{\lambda}=0.
\end{align}When $t_0\geq -2\text{Re}\; (t_1\sigma)> 0$ or
$-t_0\geq 2\text{Re}\; (t_1\sigma)> 0$, this gives
\begin{align*}
\lambda=\log
\left(-\frac{t_0}{t_1\sigma+\bar{t}_1\bar{\sigma}}\right)\geq 0.
\end{align*}Substituting into the first equation in \eqref{so2} and
the first equation in \eqref{so3}, we find that the power series
\begin{align*}
&\left.\mL(w;\boldsymbol{t})\right|_{t_n=0 \, (|n|\geq
2)}\\=&-\sigma-\frac{t_0w}{4\text{Re}\;(t_1\sigma)}
\left(\left(1+\frac{\sigma}{w}\right)^2
+\left(1+\frac{\sigma}{w}\right)
\sqrt{1+\frac{2\sigma(t_0+4\text{Re}\;(t_1\sigma))}{wt_0}+\frac{\sigma^2}{w^{2}}}\right)\\
&\left.\tilde{\mL}(w;\boldsymbol{t})\right|_{t_n=0 \, (|n|\geq
2)}\\=&-\sigma-\frac{t_0}{4\text{Re}\;(t_1\sigma)w}\left(\left(\sigma+w\right)^2
-\sigma\left(\sigma+w\right)
\sqrt{1+\frac{2\bar{\sigma}(t_0+4\text{Re}\;(t_1\sigma))w}{t_0}+\bar{\sigma}^2w^{2}}\right)\\
\end{align*}satisfy the dToda hierarchy \eqref{Lax1} with
$t_{-n}=-\bar{t}_n$ for $n\geq 1$.

When $t_{-n}=-\bar{t}_n$ $\forall n\geq 1$ and $t_n=0$ for
$|n|\geq 3$, the hodograph relation \eqref{ho3} with
$R(\lambda)=0$ reads as
\begin{align*}
t_0+t_1 \sigma
e^{\lambda}+\bar{t}_1\bar{\sigma}e^{\lambda}+
2t_2\sigma^2e^{\lambda}(3e^{\lambda}-2)+2\bar{t}_2
\bar{\sigma}^2e^{\lambda}(3e^{\lambda}-2)=0,
\end{align*}
which can be solved for
$\left.\lambda(\boldsymbol{t})\right|_{t_n=0 (|n|\geq 3)}$ in a
certain domain of $\boldsymbol{t}$.
\begin{remark}
If we do not impose the condition $|\sigma|=1$ in the example
above, then the maps $G_{\lambda}(w)$ and $F_{\lambda}(w)$ still
solve the L\"owner equations \eqref{radialout}, \eqref{radialin}
with initial conditions $G_0(w)=w$ and $F_0(w)=w$ respectively. In
this case, $G_{\lambda}(w)$ maps the disc $\{ |w|\geq |\sigma|\}$
conformally onto $\{ |z|\geq |\sigma|\}\setminus \{
\alpha\sigma\;:\;\alpha\in[1,
2e^{\lambda}-1+2\sqrt{e^{2\lambda}-e^{\lambda}}]\}$ and
$F_{\lambda}(w)$ maps the disc $\{ |w|\leq |\sigma|\}$ conformally
onto $\{ |z|\leq |\sigma|\}\setminus\{ \alpha\sigma\;:\;\alpha\in[
2e^{\lambda}-1-2\sqrt{e^{2\lambda}-e^{\lambda}},1]\}.$ However,
\eqref{inver} no longer holds. Nevertheless, one can still show as
in Proposition \ref{Pr2} that with $\lambda(\boldsymbol{t})$
satisfying \eqref{hodograph2}, the conclusion of Proposition
\ref{Pr2} still holds.
\end{remark}
\subsubsection{} \textbf{Example 2} \\
Given a point $\sigma\in S^1$, let
\begin{align*}
G(w,\lambda)= e^{\lambda}\left( w+ 2\sigma +\sigma^2
w^{-1}\right)=\frac{e^{\lambda}(w+\sigma)^2}{w}.
\end{align*}It is easy to verify directly that $G(w,\lambda)$
satisfies the radial L\"owner equation \eqref{radialout} with
$\sigma(\lambda)=\sigma$ and $\phi(\lambda)=\lambda$. Therefore,
$G(w,\lambda)$ here satisfies the same equation as the
$G(w,\lambda)$ in \textbf{Example 1} but with a different initial
condition. It maps the exterior disc conformally onto
$$\hat{\C}\setminus\{ \sigma\alpha\;:\;\alpha\in [0, 4e^{\lambda}].\}
$$The function
\begin{align*}
F(w,\lambda)=\ov{G(1/\bar{w},\lambda)}^{\;-1}=
\frac{e^{-\lambda}w}{(1+\bar{\sigma}w)^2}=
e^{-\lambda}\left(w+2\bar{\sigma}w^2 +
3\bar{\sigma}^2w^3+\ldots\right)
\end{align*}maps $\Del$ conformally onto
$$\C\setminus \{ \sigma \alpha\;:\; \alpha\in [e^{-\lambda}/4,
\infty)\}.$$For $n\geq 1$, let $\Phi_n(w,\lambda)$ and
$\Psi_n(w,\lambda)$ be the Faber polynomials of the inverse
functions $g(z,\lambda), f(z,\lambda)$ of $G(w,\lambda),
F(w,\lambda)$ respectively. Then
\begin{align*}
\Phi_1(w,\lambda) =& e^{\lambda}( w + 2\sigma ), \hspace{3cm}
\Psi_1(w,\lambda)= e^{\lambda}(w^{-1} + 2\bar{\sigma}),\\
\Phi_2(w,\lambda)=& e^{2\lambda}(w^2+4\sigma w+ 6\sigma^2),
\hspace{1cm} \Psi_2(w,\lambda)=e^{2\lambda}(w^2+4\bar{\sigma} w+
6\bar{\sigma}^2).
\end{align*}
Therefore, the functions $\xi_n(\lambda)$, $n=\pm 1, \pm 2$ in
\eqref{hodograph2} are given by
\begin{align*}
\xi_1(\lambda)&= e^{\lambda}\sigma, \hspace{2cm}
\xi_{-1}(\lambda)= -e^{\lambda}\bar{\sigma},
\\
\xi_{2}(\lambda)&= 6e^{2\lambda}\sigma^2,
\hspace{1.3cm}\xi_{-2}(\lambda)=-6e^{2\lambda}\bar{\sigma}^2.
\end{align*}When $t_{-1}=-\bar{t}_1$ and $t_n=0$ for $|n|\geq 2$,
the hodograph relation
\eqref{ho3} with $R(\lambda)=0$ is the same as equation
\eqref{sol1}. Therefore, when $t_0\geq -2\text{Re}\; (t_1\sigma)>0
$ or $-t_0\geq 2\text{Re}\; (t_1\sigma)> 0$, we find that the
power series
\begin{align*}
\left.\mL(w;\boldsymbol{t})\right|_{t_n=0 \, (|n|\geq
2)}=&-\frac{t_0}{2\text{Re}\; (t_1\sigma)}
\left(w+2\sigma+\sigma^2w^{-1}\right)\\
\left.\tilde{\mL}(w;\boldsymbol{t})\right|_{t_n=0 \, (|n|\geq
2)}=&-\frac{2\text{Re}\;
(t_1\sigma)}{t_0}\frac{w}{(1+\bar{\sigma}w)^2}
\end{align*}satisfies the dToda hierarchy \eqref{Lax1} with
$t_{-n}=-\bar{t}_n$ for $n\geq 1$.
When $t_{-n}=-\bar{t}_n$
$\forall n\geq 1$ and $t_n=0$ for $|n|\geq 3$, the hodograph
relation \eqref{ho3} with $R(\lambda)=0$ reads as
\begin{align*}
t_0+t_1 \sigma e^{\lambda}+\bar{t}_1\bar{\sigma}e^{\lambda}+
6t_2e^{2\lambda}\sigma^2+ 6\bar{t}_2e^{2\lambda}\bar{\sigma}^2=0,
\end{align*} which gives
\begin{align*}
\left.\lambda(\boldsymbol{t})\right|_{t_n=0 (|n|\geq 3)} =\log
\left( \frac{-\text{Re}\; (t_1\sigma) +\sqrt{[\text{Re}\;
(t_1\sigma)]^2 - 12t_0\text{Re}\; (t_2\sigma^2)}}{12\text{Re}\;
(t_2\sigma^2)}\right)
\end{align*}in a certain domain of $\boldsymbol{t}$.


\begin{thebibliography}{10}

\bibitem{BB} M.~Bauer and D.~Bernard,
\emph{2D growth processes: SLE and Loewner
chains}, e-print archive: math-ph/0602049.

\bibitem{Bie}
L.~Bieberbach, \emph{\"{U}ber die {K}oeffizienten derjenigen
{P}otenzreihen,
  welche eine schlichte {A}bbildung des {E}inheitskreises vermitteln},
  S.-B.Preuss.Akad.Wiss. (1916), 940--955.

\bibitem{ref1}
A.~Boyarsky, A.~Marshakov, O.~Ruchayskiy, P.~Wiegmann, and
A.~Zabrodin,
  \emph{Associativity equations in dispersionless integrable hierarchies},
  Phys. Lett. B \textbf{515} (2001), no.~3-4, 483--492.

\bibitem{Cardy} J.~Cardy, \emph{SLE for theoretical physicists},
Annals Phys. \textbf{318} (2005) 81--118.

\bibitem{ref13_2}
R.~Carroll and Y.~Kodama, \emph{Solution of the dispersionless
{H}irota
  equations}, J. Phys. A \textbf{28} (1995), no.~22, 6373--6387.

\bibitem{CY}
Ling-Lie Chau and Yue Yu, \emph{Unitary polynomials in normal
matrix models and
  wave functions for the fractional quantum {H}all effects}, Phys. Lett. A
  \textbf{167} (1992), no.~5-6, 452--458.

\bibitem{ref20}
Ling-Lie Chau and Oleg Zaboronsky, \emph{On the structure of
correlation
  functions in the normal matrix model}, Comm. Math. Phys.
\textbf{196} (1998),
  no.~1, 203--247.

\bibitem{Conway}
John~B. Conway, \emph{Functions of one complex variable. {II}.},
Graduate Texts
  in Mathematics, 159. Springer-Verlag, New York, 1995.

\bibitem{deB}
L.~De~Branges, \emph{A proof of the {B}ieberbach conjecture}, Acta
Math.
  \textbf{154} (1985), no.~1-2, 137--152.

\bibitem{ref15}
B.~A. Dubrovin and S.~M. Natanzon, \emph{Real theta-function
solutions of the
  {K}adomtsev-{P}etviashvili equation}, Izv. Akad. Nauk SSSR Ser. Mat.
  \textbf{52} (1988), no.~2, 267--286, 446.

\bibitem{Duren}
P.~L. Duren, \emph{Univalent functions}, Grundlehren der
Mathematischen
  Wissenschaften [Fundamental Principles of Mathematical Sciences], vol. 259,
  Springer-Verlag, New York, 1983.

\bibitem{ref13_1}
John Gibbons and Yuji Kodama, \emph{Solving dispersionless {L}ax
equations},
  Singular limits of dispersive waves (Lyon, 1991), NATO Adv. Sci. Inst. Ser. B
  Phys., vol. 320, Plenum, New York, 1994, pp.~61--66.

\bibitem{ref10}
John Gibbons and Serguei~P. Tsarev, \emph{Conformal maps and
reductions of the
  {B}enney equations}, Phys. Lett. A \textbf{258} (1999), no.~4-6, 263--271.

\bibitem{ref19}
S.~Kharchev, A.~Marshakov, A.~Mironov, and A.~Morozov,
\emph{Generalized
  {K}ontsevich model versus {T}oda hierarchy and discrete matrix models},
  Nuclear Phys. B \textbf{397} (1993), no.~1-2, 339--378.

\bibitem{ref7}
I.~K. Kostov, I.~Krichever, M.~Mineev-Weinstein, P.~B. Wiegmann,
and
  A.~Zabrodin, \emph{The {$\tau$}-function for analytic curves}, Random matrix
  models and their applications, Math. Sci. Res. Inst. Publ., vol.~40,
  Cambridge Univ. Press, Cambridge, 2001, pp.~285--299.

\bibitem{Kufarev}
P.~P. Kufarev, \emph{A remark on integrals of {L}\"owner's
equation}, Doklady
  Akad. Nauk SSSR (N.S.) \textbf{57} (1947), 655--656.

\bibitem{LSW}
Gregory~F. Lawler, Oded Schramm, and Wendelin Werner, \emph{Values
of
  {B}rownian intersection exponents. {I}. {H}alf-plane exponents}, Acta Math.
  \textbf{187} (2001), no.~2, 237--273.

\bibitem{Lo}
K.~L\"owner, \emph{{U}ntersuchungen uber schlichte konforme
{A}bbildungen des
  {E}inheitskreises. i.}, Math. Ann. \textbf{89} (1923), 103--121.

\bibitem{Manas}
Manuel Ma{\~n}as, \emph{{$S$}-functions, reductions and hodograph
solutions of
  the {$r$}th dispersionless modified {KP} and {D}ym hierarchies}, J. Phys. A
  \textbf{37} (2004), no.~46, 11191--11221.

\bibitem{MMAM}
Manuel Ma{\~n}as, Luis Mart{\'{\i}}nez~Alonso, and Elena Medina,
  \emph{Reductions and hodograph solutions of the dispersionless {KP}
  hierarchy}, J. Phys. A \textbf{35} (2002), no.~2, 401--417.

\bibitem{ref8}
A.~Marshakov, P.~Wiegmann, and A.~Zabrodin, \emph{Integrable
structure of the
  {D}irichlet boundary problem in two dimensions}, Comm. Math. Phys.
  \textbf{227} (2002), no.~1, 131--153.

\bibitem{MarshallRohde}
 Donald E.~Marshall and Steffen Rohde,
\emph{The Loewner differential equation and slit mappings}.
J. Amer. Math. Soc. \textbf{18} (2005), no. 4, 763--778


\bibitem{ref6_1}
M.~Mineev-Weinstein, P.B. Wiegmann, and A.~Zabrodin,
\emph{Integrable sructure
  of interface dynamics}, Phys.Rev.Lett. \textbf{84} (2000), 5106--5109.

\bibitem{Pom}
C.~Pommerenke, \emph{Univalent functions}, Vandenhoeck \&
Ruprecht,
  G\"ottingen, 1975, With a chapter on quadratic differentials by Gerd Jensen,
  Studia Mathematica/Mathematische Lehrb\"ucher, Band XXV.

\bibitem{Schramm}
Oded Schramm, \emph{Scaling limits of loop-erased random walks and
uniform
  spanning trees}, Israel J. Math. \textbf{118} (2000), 221--288.

\bibitem{TT2}
K.~Takasaki and T.~Takebe, \emph{{${\rm SDiff}(2)$} {KP}
hierarchy}, Infinite
  analysis, Part A, B (Kyoto, 1991), Adv. Ser. Math. Phys., vol.~16, World Sci.
  Publishing, River Edge, NJ, 1992, pp.~889--922.

\bibitem{TT4}
\bysame, \emph{Quasi-classical limit of {T}oda hierarchy and
{$W$}-infinity
  symmetries}, Lett. Math. Phys. \textbf{28} (1993), no.~3, 165--176.

\bibitem{TT1}
\bysame, \emph{Integrable hierarchies and dispersionless limit},
Rev. Math.
  Phys. \textbf{7} (1995), no.~5, 743--808.

\bibitem{TTnew}
\bysame, \emph{Radial {L}\"owner equation and dispersionless
cm{KP} hierarchy},
  arXiv:nlin.SI/0601063 (2006).

\bibitem{dcmKP} Lee-Peng Teo, \emph{On dispersionless coupled
modified KP hierarchy}, e-print archive: nlin.SI/0304007.

\bibitem{Teo}
Lee-Peng Teo, \emph{Analytic functions and integrable
  hierarchies---characterization of tau functions}, Lett. Math. Phys.
  \textbf{64} (2003), no.~1, 75--92.

\bibitem{Tsarev}
S.~P. Tsar{\"e}v, \emph{The geometry of {H}amiltonian systems of
hydrodynamic
  type. {T}he generalized hodograph method}, Izv. Akad. Nauk SSSR Ser. Mat.
  \textbf{54} (1990), no.~5, 1048--1068.

\bibitem{Wein}
Lenard Weinstein, \emph{The {B}ieberbach conjecture}, Internat.
Math. Res.
  Notices (1991), no.~5, 61--64.

\bibitem{ref6_2}
P.~B. Wiegmann and A.~Zabrodin, \emph{Conformal maps and
integrable
  hierarchies}, Comm. Math. Phys. \textbf{213} (2000), no.~3, 523--538.

\bibitem{ref11}
L.~Yu and J.~Gibbons, \emph{The initial value problem for
reductions of the
  {B}enney equations}, Inverse Problems \textbf{16} (2000), no.~3, 605--618.

\bibitem{ref9}
A.~V. Zabrodin, \emph{The dispersionless limit of the {H}irota
equations in
  some problems of complex analysis}, Teoret. Mat. Fiz. \textbf{129} (2001),
  no.~2, 239--257.

\end{thebibliography}

\end{document}